\documentclass[a4paper]{article}
\title{Inequalities for the Lattice Width of Lattice-Free\\ Convex Sets in the Plane} 
\author{Gennadiy Averkov and Christian Wagner}

\usepackage{amsfonts,amsthm,amsmath,amssymb,enumerate}
\usepackage{graphicx}
\usepackage{graphpap}
\usepackage{color}
\usepackage{subfigure}
\usepackage{empheq}
\usepackage[bf]{caption}
\usepackage{srcltx}

\numberwithin{equation}{section}
\numberwithin{subfigure}{section}

\newcommand{\aff}{\rmcmd{aff}}
\newcommand{\area}{A}
\newcommand{\bd}{\rmcmd{bd}}
\newcommand{\ceil}[1]{\left\lceil#1\right\rceil}
\newcommand{\conv}{\rmcmd{conv}}
\newcommand{\dotvar}{\,\cdot\,}
\newcommand{\floor}[1]{\left\lfloor #1 \right\rfloor}

\newcommand{\integer}{\mathbb{Z}}
\newcommand{\rational}{\mathbb{Q}}
\newcommand{\intr}{\rmcmd{int}}
\newcommand{\KK}{\mathcal{K}}
\newcommand{\modulo}[1]{\ (\mathrm{mod} \ #1)}

\newcommand{\natur}{\mathbb{N}}
\newcommand{\nint}[1]{\left\lfloor #1 \right\rceil}

\newcommand{\overtwocond}[2]{\substack{{#1} \\ {#2}}}

\newcommand{\real}{\mathbb{R}}

\newcommand{\rmcmd}[1]{\mathop{\mathrm{#1}}}
\newcommand{\setcond}[2]{\left\{#1\,:\,#2\right\}}

\newcommand{\term}[1]{\emph{#1}}
\newcommand{\thmcitation}[1]{{\upshape #1}}

\newcommand{\lwid}{w}
\newcommand{\wid}{w}

\newtheorem{nn}{}[section]
\newtheorem{corollary}[nn]{Corollary}
\newtheorem{lemma}[nn]{Lemma}
\newtheorem*{lemma*}{Lemma}

\newtheorem{theorem}[nn]{Theorem}
\newtheorem{proposition}[nn]{Proposition}

\begin{document}
\maketitle

\begin{abstract}
  A closed, convex set $K$ in $\real^2$ with non-empty interior is called lattice-free if the interior of $K$ is disjoint with $\integer^2$. In this paper we study the relation between the area and the
  lattice width of a planar lattice-free convex set in the general and
  centrally symmetric case.
  A correspondence between lattice width on the one hand
  and covering minima on the other, allows us to reformulate our
  results in terms of covering minima introduced by Kannan and
  Lov\'asz \cite{KannanLovasz88}.
  We obtain a sharp upper bound for the area for any given value of
  the lattice width.
  The lattice-free convex sets satisfying the upper bound are
  characterized.
  Lower bounds are studied as well.
  Parts of our results are applied in \cite{AWW09} for cutting plane
  generation in mixed integer linear optimization, which was the original
  inducement for this paper.
  %for a classification of three-dimensional maximal
  %lattice-free polytopes with integer vertices.
  %These polytopes are used in (mixed) integer linear optimization
  %for deriving cutting planes.
  We further rectify a result of \cite{KannanLovasz88} with a new
  proof.
\end{abstract}

\newtheoremstyle{itsemicolon}{}{}{\mdseries\rmfamily}{}{\itshape}{:}{ }{}
\theoremstyle{itsemicolon}
\newtheorem*{msc*}{Mathematics Subject Classification (AMS 2010)}
\begin{msc*}
        Primary: 52C05; Secondary: 52A38, 52A40, 52C15
\end{msc*}

\newtheorem*{keywords*}{Key words and phrases}

\begin{keywords*}
  area; convex set; covering minimum; inhomogeneous minimum; lattice-free
  body; lattice width
\end{keywords*}

%%%%%%%%%%%%%%%%%%%%%%%%%%%%%%%%%%%%%%%%%%%%%%%%%%%%%%%%%%%%%%%%%%%%%%
%%%%%%%%%%%%%%%%%%%%%%%%%%%%%%%%%%%%%%%%%%%%%%%%%%%%%%%%%%%%%%%%%%%%%%

\section{Introduction} \label{intro}

% IDEA: a ``user-friendly'' introduction to lattice width and covering minima

This paper is devoted to the relation between the area and the lattice
width (resp.~the area and the covering minima) of a lattice-free
two-dimensional convex set.
Let us denote by $\KK^2$ the class of closed, convex sets in $\real^2$
with non-empty interior.
We call a set $K \in \KK^2$ \term{lattice-free} if the interior of $K$
is disjoint with $\integer^2$.
A convex set $S$ is said to be a \term{strip} if $S=\conv (l^1 \cup
l^2)$, where $l^1$ and $l^2$ are distinct parallel lines.
The Euclidean distance between $l^1$ and $l^2$ is said to be the width
of $S$.
If the interior of $S$ does not contain integer points and the sets
$l^1 \cap \integer^2$, $l^2 \cap \integer^2$ are affine images of
$\integer$, then we call $S$ a \term{split} since $S$ splits
$\integer^2$ into two parts.
Let $K \in \KK^2$. If $S=\conv (l^1\cup l^2)$ is a split and $l_K^1$ and $l_K^2$  are distinct supporting lines of $K$ parallel to $l^1$ and
$l^2$, we define the strip $S_K:=\conv (l^1_K \cup l^2_K)$.
Then the \term{lattice width} $\lwid(K)$ of $K$ is defined as the minimum of the ratio
between the width of $S_K$  and the width of $S$ among all splits $S$, for which $S_K$ exists.
In analytic terms, $\lwid(K) := \min \setcond{\lwid(K,u)}{u \in \integer^2 \setminus
  \{o\}} $, where $\lwid(K,u)$ is the \term{width function}
defined by $\lwid(K,u):=\max\{u^\top x : x \in K\} - \min\{u^\top x :
x \in K\}$ for $u \in \real^2$.

With these notions our main contribution is the following:
for a given lattice-free set $K \in \KK^2$ we present a list of
inequalities which relate its area $A(K)$ to its lattice width
$\wid(K)$, see Theorems~\ref{AreaLWidth} and
\ref{CentrSym} below. 
The inequalities which give upper bounds for $A(K)$ for a given
$\lwid(K)$, and the sets yielding equality in these inequalities are
characterized.
For the case of centrally symmetric sets $K \in \KK^2$ we even give
the complete list of inequalities, that is, the lower and upper bounds
for $A(K)$ for a given $\lwid(K)$ and a characterization of all pairs
$(\lwid(K),A(K))$ where equality is attained.

All results obtained in this paper can be formulated in terms of an
arbitrary lattice, but for the sake of simplicity we use
$\integer^2$.
Let $\Lambda$ be an arbitrary lattice in $\real^2$.
We can introduce the lattice width $\lwid(K,\Lambda)$ of $K$
with respect to $\Lambda$ (for the precise definition see
\cite{KannanLovasz86,KannanLovasz88}).
Then $\lwid(T(K),T(\Lambda))=\lwid(K,\Lambda)$ for every
linear transformation $T$ in $\real^2$.
Choosing $T$ such that it maps $\Lambda$ onto $\integer^2$
we obtain $\lwid(T(K)) = \lwid(T(K),\integer^2) =
\lwid(T(K),T(\Lambda)) = \lwid(K,\Lambda)$.
Since $A(T(K)) \cdot \det \Lambda = A(K)$, every inequality
involving the area and the lattice width with respect to
$\integer^2$ can be transformed to an inequality involving
the area, the lattice width  with respect to $\Lambda$ and
$\det \Lambda$.
For figures it is sometimes more convenient to use the
lattice of regular triangles, i.e., the lattice generated by
the vectors $(1,0)$ and $(\frac{1}{2},\frac{\sqrt{3}}{2})$.
For information on lattices and convexity we refer to
\cite{GruberBook07} and \cite{GruLek87}.

Our motivation for studying these relations was the application of
our results for a classification of three-dimensional lattice-free
polyhedra having an integer point in the relative interior of each
facet, see \cite{AWW09}. 
Such class of polyhedra can be used in mixed integer linear
programming for deriving cutting planes; see also
\cite[Chapter\,23]{schrijver-book-86} for background information on
cutting plane theory.

The \term{second covering minimum} $\mu_2(K)$ of $K \in \KK^2$ is
defined as the minimal $t \ge 0$ such that the sets $t K +
\integer^2$ cover $\real^2$. The value $\mu_2(K)$ is also
known under the name of \term{inhomogeneous minimum}, see
\cite[p.\,98]{GruLek87}. It turns out that some translate of
$K$ is lattice-free if and only if $\mu_2(K) \ge 1$. The
\term{first covering minimum} $\mu_1(K)$ is defined as the
minimal $t \ge 0$ such that every line in $\real^2$
intersects $t K + \integer^2$. One can show that $\mu_1(K)
\lwid(K) = 1$ for every $K \in \KK^2$, see e.g.~\cite[Lemma\,2.3]{KannanLovasz88}. This leads to a correspondence between the lattice width on the one hand and the covering minima on the other, provided that $K$ is lattice-free. The results we present in this paper can therefore be expressed as a
relation between the area and the covering minima of $K$, as well; see
Corollaries~\ref{Upper bounds for area in terms of mu's} and
\ref{mu.sym}.
The notions lattice width and covering minima were introduced by Kannan and
Lov\'asz \cite{KannanLovasz86,KannanLovasz88} for an arbitrary dimension; see also the papers of Khinchin \cite{Khinchin48} and Fejes
T\'oth and Makai Jr.~\cite{MakaiToth74} for the earlier related results.

The paper has the following structure. Our main results
are stated in Section \ref{results}. 
Basic notions and results which we shall need to prove our theorems are
introduced in Section \ref{prelims}.
Section~\ref{aux:tri} presents formulas for the lattice width and area
of triangles.
Section \ref{general.case} contains the proofs for general planar
lattice-free convex sets.
The proofs for the centrally symmetric planar lattice-free convex sets
are given in Section \ref{centr.sym.case}.

%%%%%%%%%%%%%%%%%%%%%%%%%%%%%%%%%%%%%%%%%%%%%%%%%%%%%%%%%%%%%%%%%%%%%%
%%%%%%%%%%%%%%%%%%%%%%%%%%%%%%%%%%%%%%%%%%%%%%%%%%%%%%%%%%%%%%%%%%%%%%

\section{Results} \label{results}

Throughout the paper, sequences with $n$ elements are indexed modulo $n$. Affine transformations preserving $\integer^2$ will be called
\term{unimodular}.

For the plane it is easy to construct examples of lattice-free convex
sets with a lattice width of two; take for example $\conv \{(0,0),
(2,0), (0,2)\}$.
What may be surprising, this value can be exceeded.
This was noticed by Hurkens \cite{Hurkens90}, who also computed the
sharp upper bound for the lattice width of a lattice-free convex set.

\begin{figure}[hbt]
  \begin{center}
  \unitlength=1mm
%       \fbox{
  \begin{picture}(40,40)
    \put(0,0){\includegraphics[width=40\unitlength]{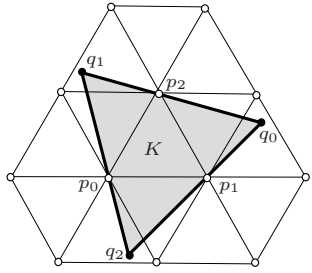}}
                \put(18,15){\scriptsize $K$}
                \put(33.2,17){\scriptsize $q_0$}
                \put(10.5,27){\scriptsize $q_1$}
                \put(13,1.5){\scriptsize $q_2$}
                \put(21,24){\scriptsize $p_2$}
                \put(9.5,10.5){\scriptsize $p_0$}
                \put(28,10.3){\scriptsize $p_1$}
%               \graphpaper[5](0,0)(40,40)
    \end{picture}
%       }
  \end{center}
  \caption{\label{hurkens-fig} Lattice-free triangle (shaded) with
    lattice width $1 + \frac{2}{\sqrt{3}}$}
\end{figure}

\begin{theorem} \thmcitation{\cite[p.\,122]{Hurkens90}}
  \label{hurkens-thm} Let $K \in \KK^2$ be lattice-free.
  Then
  \begin{equation}\label{LWidthUpper}
    \lwid(K) \le 1+ \frac{2}{\sqrt{3}},
  \end{equation}
  with equality if and only if $K$ is a triangle with vertices
  $q_0,q_1,q_2$ such that, for every $i$,  the point
  $p_i:=\frac{1}{\sqrt{3}} q_{i+1} +\left(1-\frac{1}{\sqrt{3}}\right) q_{i+2}$
  belongs to $\integer^2$ (see Fig.~\ref{hurkens-fig}).
\end{theorem}

Theorem \ref{AreaLWidth} states the relation between the area and
the lattice width of arbitrary lattice-free convex sets in the plane.

\begin{figure}[ht]
  \unitlength=1mm
  \centering
  \subfigure[\label{quad-area-max-fig} Illustration to
  \eqref{AreaEasy} and Part \ref{AreaEasyEq}]{
                \begin{picture}(30,35)
                        \put(0,0){\includegraphics[width=30mm]{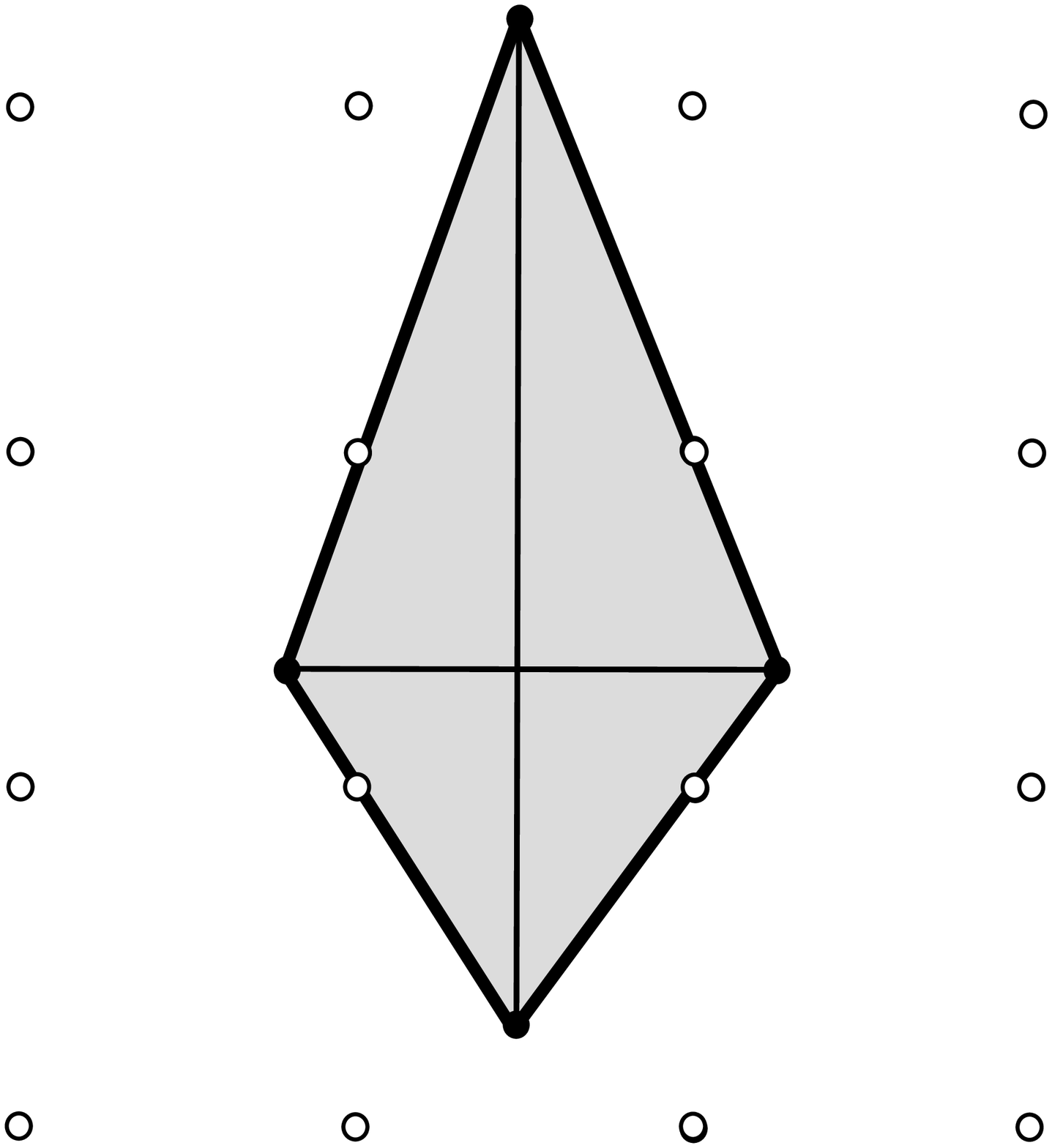}}
                        \put(17,14.5){\scriptsize $I_1$}
                        \put(15,20){\scriptsize $I_2$}
                        \put(11,16){\scriptsize $K$}
                \end{picture}
        } 
  \qquad 
  \subfigure[\label{tri-area-max-fig} Illustration to
  \eqref{AreaHard} and Part \ref{AreaHardEq}]{
        \begin{picture}(35,35)
        \put(2,4){\includegraphics[width=32\unitlength]{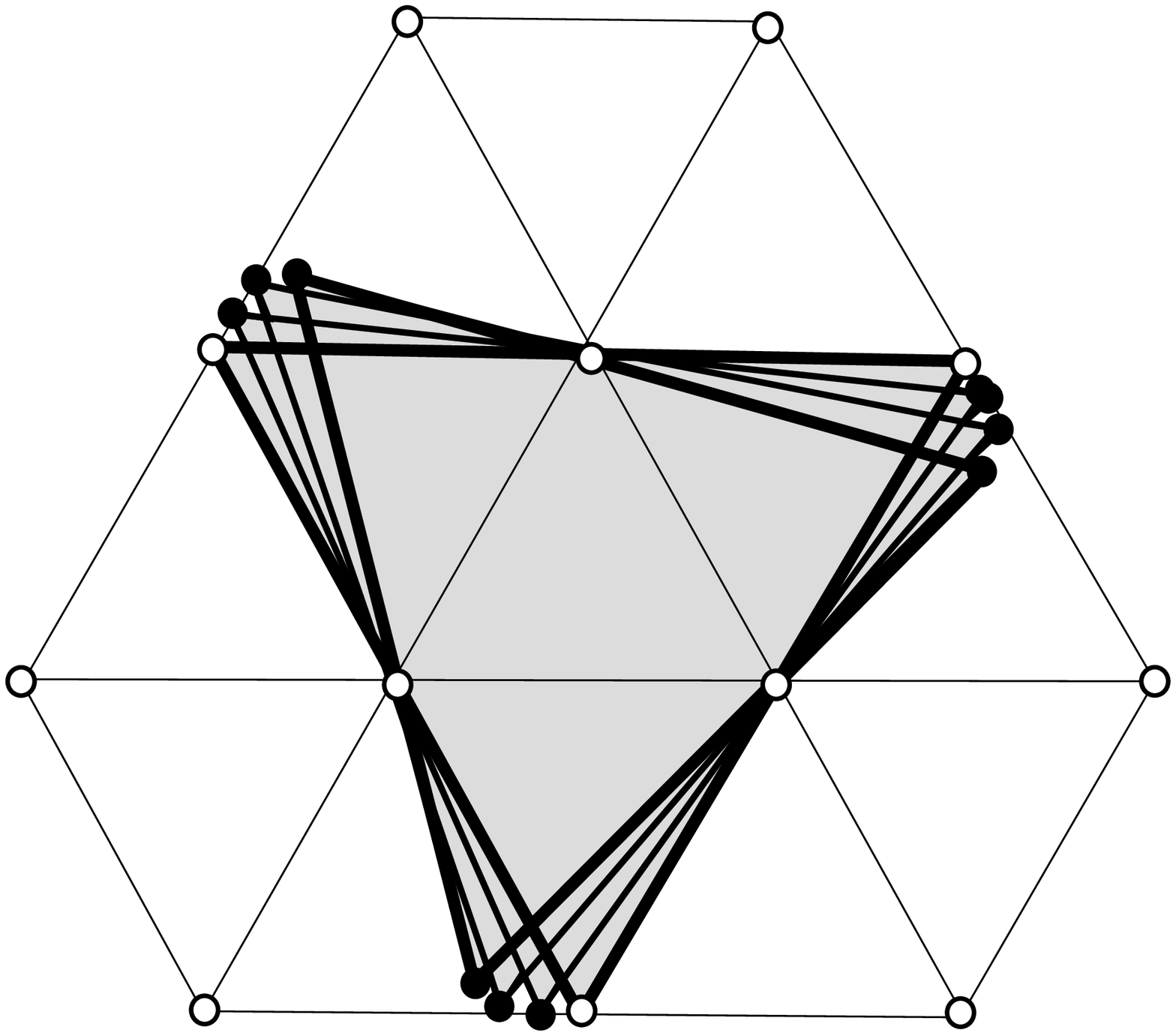}}
                        \put(16.8,16){\scriptsize $K$}
                        \put(10.5,25.5){\scriptsize $q_1$}
                        \put(12,5.5){\scriptsize $q_2$}
                        \put(28,17){\scriptsize $q_0$}
                        \put(9,14.5){\scriptsize $p_0$}
                        \put(24.7,12){\scriptsize $p_1$}
                        \put(19.2,23.6){\scriptsize $p_2$}
        \end{picture}
        } 
  \qquad 
  \subfigure[\label{area-min-fig} Illustration to
    \eqref{AreaBelow} and Part \ref{AreaBelowEq}]{
               \begin{picture}(60,35)
                                                                \put(1,2){\includegraphics[width=57\unitlength]{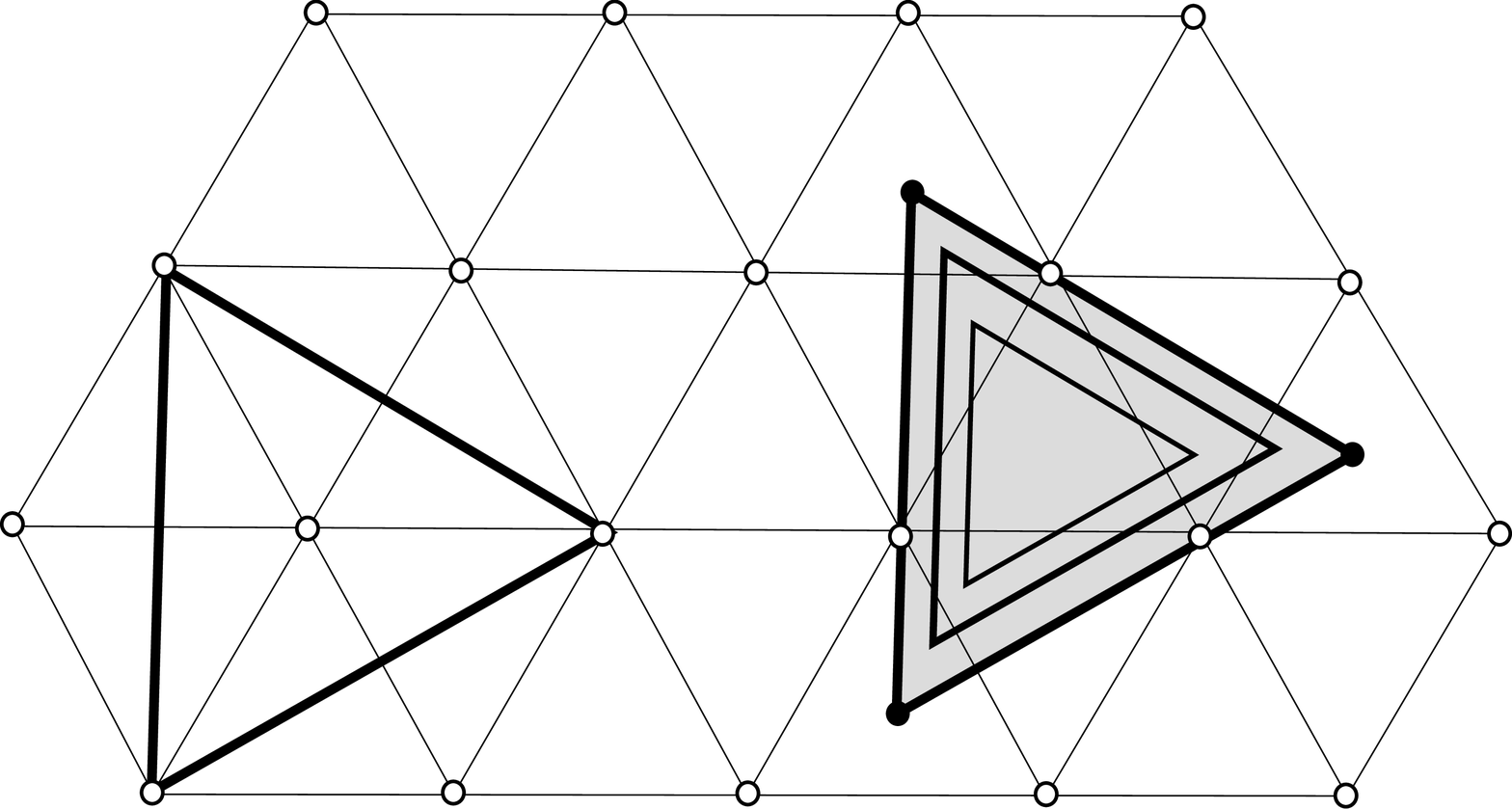}}
                        \put(0,0){\tiny $(-1,-1)$}
                        \put(25,10.6){\tiny $(1,0)$}
                        \put(1,22){\tiny $(0,1)$}
                        \put(14,10.6){\tiny $(0,0)$}
                        \put(39.2,14.7){\tiny $K$}
               \end{picture}
        }
  \parbox[t]{0.80\textwidth}{\caption{Examples of sets yielding equality in
    \eqref{AreaEasy}--\eqref{AreaBelow} (shaded) \label{optimal-bodies-fig}}}
\end{figure}

\begin{figure}[ht]
        \unitlength=1mm
   \centering
   \begin{picture}(40,32)
      \put(2,2){\includegraphics[width=41\unitlength]{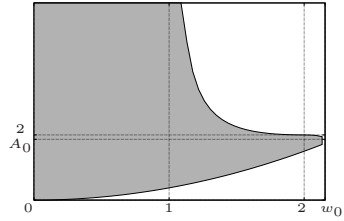}}
      \put(2,2){\tiny $0$}
      \put(20.5,2){\tiny $1$}
      \put(38,2){\tiny $2$}
      \put(41,2){\tiny $w_0$}
      \put(0,10.5){\tiny $A_0$}
      \put(1,12.5){\tiny $2$} 
   \end{picture} \\
        \parbox[t]{0.80\textwidth}{\caption{\label{diagram-general} Pairs $(w,A)$ satisfying the inequalities of Theorem~\ref{AreaLWidth}; $w_0 = 1+2/\sqrt{3}$, $A_0$ is the area of $K$ with $\lwid(K)=w_0$}}
\end{figure}

\begin{theorem} \label{AreaLWidth}
Let $K \in \KK^2$ be lattice-free with  $\wid:= \lwid(K)$ and $A :=
\area(K)$.
Then
\begin{align}
  A & \le \infty & & \mbox{for} & & 0 < \wid \le
        1,\label{InfBound} \\
  A & \le \frac{\wid^2}{2 (\wid -1 )} & &\mbox{for} & & 1 <\wid
        \le 2, \label{AreaEasy} \\
  A & \le \frac{3 \wid^2}{3 \wid + 1 - \sqrt{1 + 6 \wid  - 3
        \lwid^2}} & &\mbox{for} & & 2  < \lwid \le 1 +
        \frac{2}{\sqrt{3}}, \label{AreaHard} \\
  A & \ge \frac{3}{8} \lwid^2 & &\mbox{for} & & 0 < \lwid \le
        1+\frac{2}{\sqrt{3}} \label{AreaBelow} 
\end{align}
(see Fig.~\ref{diagram-general}). Furthermore, the following statements hold.
\begin{enumerate}[I.]
  \item \label{AreaInfEq} Equality in \eqref{InfBound} is
    attained if and only if $K$ is unbounded and contained in a
    split.
  \item \label{AreaEasyEq} Equality in
    \eqref{AreaEasy} is attained if and only if, up to unimodular
    transformations, $K = \conv (I_1 \cup I_2)$, where $I_1$ is a
    translate of $\conv \{(0,0),(\lwid,0) \}$, $I_2$ is a
    translate of $\conv \{(0,0),(0,\frac{\lwid}{\lwid-1})\}$, and
    $I_1 \cap I_2 \not = \emptyset$ (see
    Fig.~\ref{quad-area-max-fig}).
  \item \label{AreaHardEq} Equality in \eqref{AreaHard} is
    attained if and only if $K$ is a triangle with vertices
    $q_0,q_1,q_2$ such that, for every $i$, the point
    $p_i:=\lambda q_{i+1} + (1-\lambda) q_{i+2}$ belongs to
    $\integer^2$ for
    \[\lambda:=\frac{3 \lwid + 1 - \sqrt{1+6 \lwid - 3
      \lwid^2}}{6 \lwid }\]
    (see Fig.~\ref{tri-area-max-fig}).
  \item \label{AreaBelowEq} If $0 < \lwid \le 2$, then equality
    in \eqref{AreaBelow} is attained if and only if, up to
    unimodular transformations, $K$ is a translate of 
    $\frac{\lwid}{2} \conv \{(1,0),(0,1),(-1,-1)\}$ (see
    Fig.~\ref{area-min-fig}).
\end{enumerate}
\end{theorem}

%\begin{figure}[hbt]
%  \begin{center}
%    \includegraphics[width=70mm]{figures/fig-th1.eps}
%  \end{center}
%  \caption{Lower and upper bound for the area (ordinate)
%    dependent on the lattice width (abscissa)}
%  \label{graph-bounds}
%\end{figure}

The bound \eqref{AreaBelow} is not sharp when $2 < \lwid \le
1+\frac{2}{\sqrt{3}}$. To see why, we need a result of Fejes T\'oth
and Makai Jr.~\cite{MakaiToth74}.
%who gave a lower bound for the area of $K \in \KK^2$ in
%terms of the lattice width.

\begin{theorem} \label{makaitoth-thm} \thmcitation{\cite{MakaiToth74}}
  Let $K \in \KK^2$ with $\lwid:=\lwid(K)$ and $A:=\area(K)$.
  Then $A \ge \frac{3}{8} \lwid^2$ with equality if and only if, up to
  unimodular transformations, $K$ is a translate of $\frac{\lwid}{2}
  \conv \{(1,0),(0,1)$, $(-1,-1)\}$.
\end{theorem}

It is easy to see that $\frac{\lwid}{2} \conv \{(1,0),(0,1),(-1,-1)\}$
does not have a lattice-free translate for $\lwid > 2$.
Thus, in view of Theorem~\ref{makaitoth-thm}, \eqref{AreaBelow} is not
sharp when $2 < \lwid \le 1+\frac{2}{\sqrt{3}}$.
The problem to find the sharp lower bound in this case is still open.

A statement analogous to Theorem~\ref{AreaLWidth} can also be proved
for the class of centrally symmetric planar convex sets.
We show the following theorem.

\begin{figure}[ht]
        \unitlength=1mm
        \centering
        \subfigure[\label{area-maximizer-sym-fig} Illustration to \eqref{SymAreaUpper} and Part \ref{eq case for area bound}]{
%               \fbox{
                        \begin{picture}(30,30)
                                \put(0,0){\includegraphics[width=30mm]{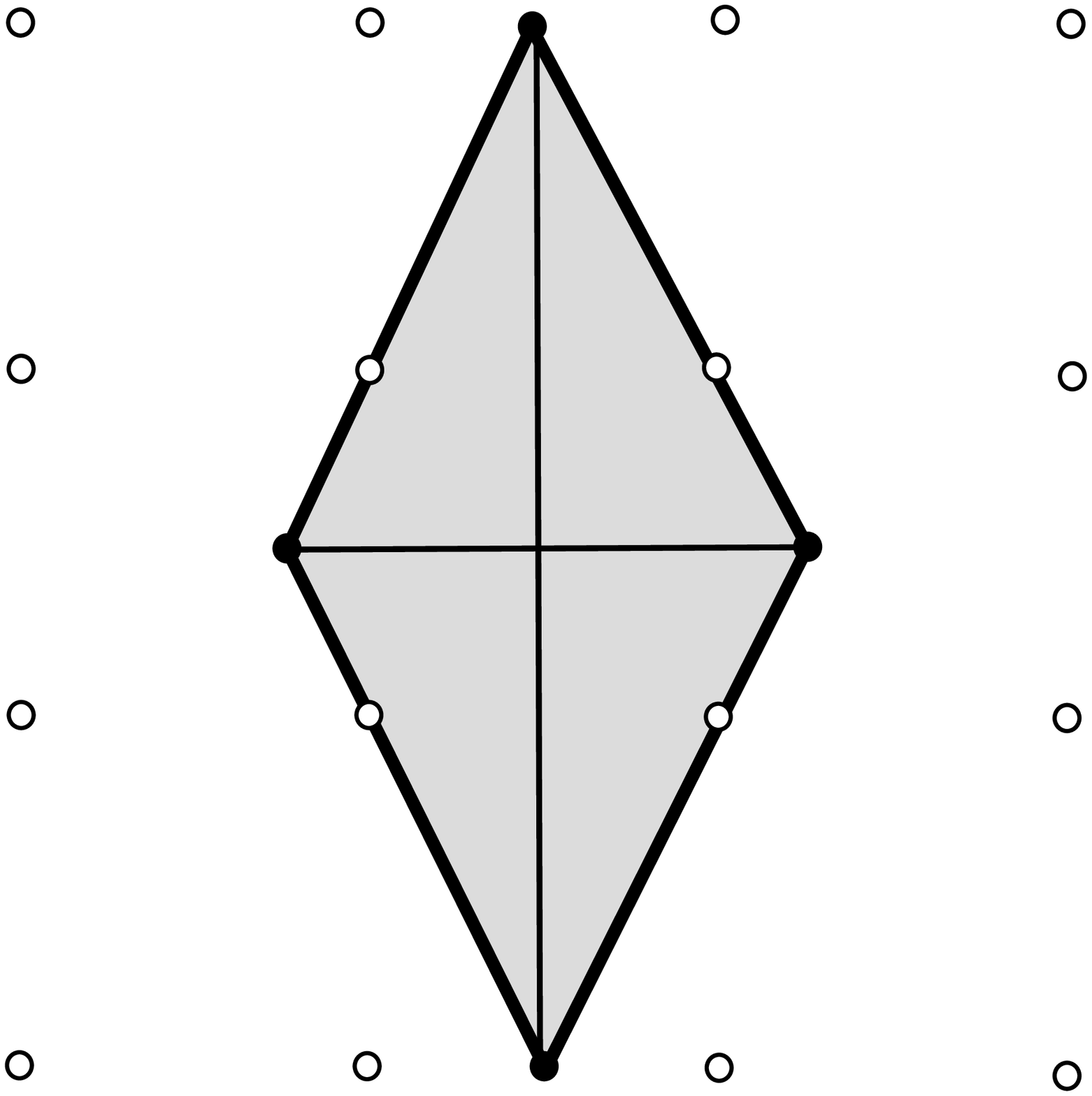}}
                                \put(16,11){\scriptsize $K$}
                                %\graphpaper[5](0,0)(30,30)
                        \end{picture}
%               }
        } 
        \qquad 
        \subfigure[\label{area-minimizer-sym-fig} Illustration to \eqref{SymAreaLower} and Part \ref{eq case lower bound sym} (with examples of $K$ for two different values of $\alpha$)]{
%               \fbox{
                        \begin{picture}(80,30)
                        \put(10,0){
                                \put(0,0){\includegraphics[width=60\unitlength]{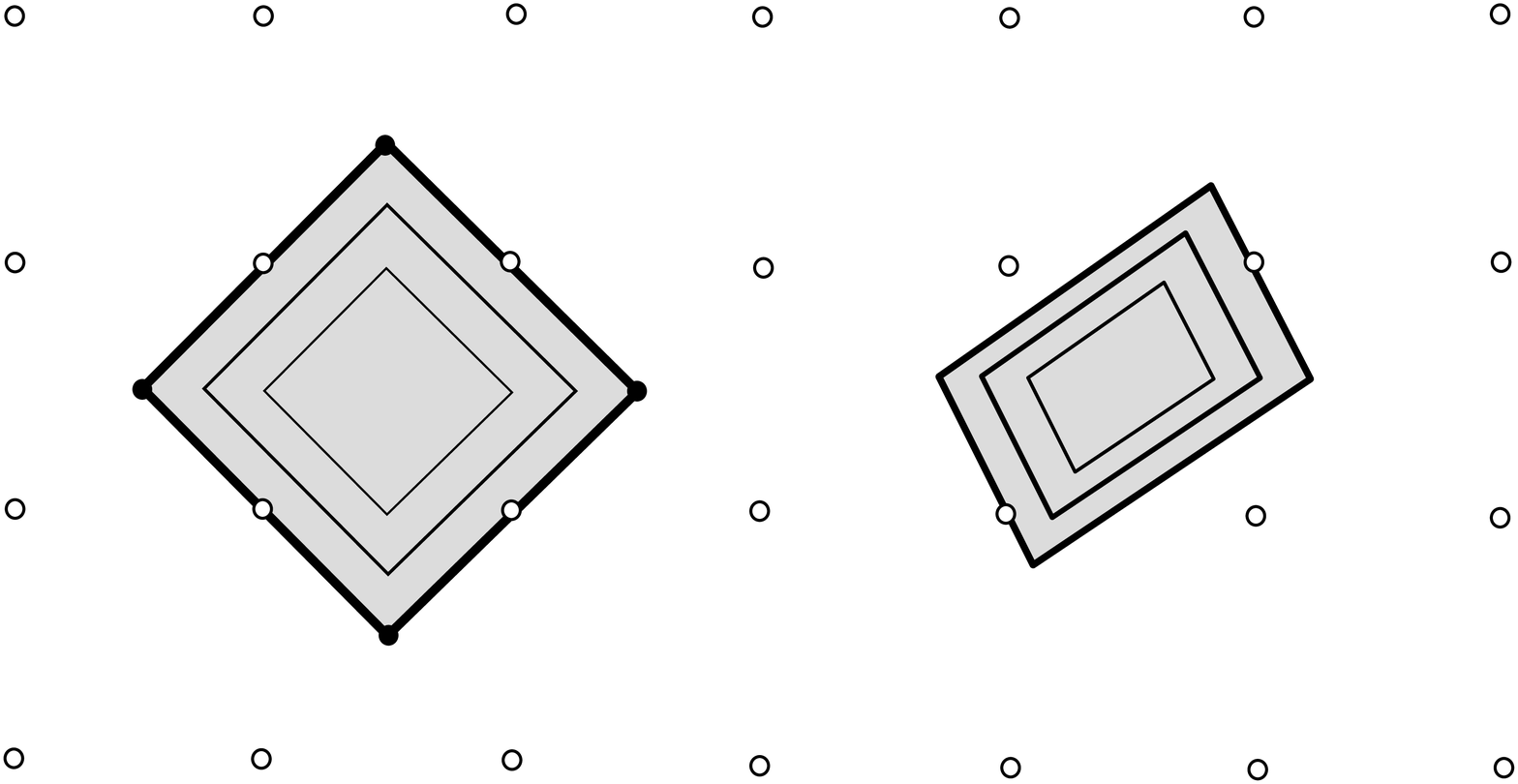}}
                                \put(14,14.5){\scriptsize $K$}
                                \put(43,15){\scriptsize $K$}
%                               \graphpaper[5](0,0)(60,30)
                        }
                        \end{picture}
%               }
        }
 \parbox[t]{0.95\textwidth}{\caption{\label{Sym.Equal} Examples of sets yielding equality in
   \eqref{SymAreaUpper} and \eqref{SymAreaLower} (shaded) }}
\end{figure}

\begin{figure}[ht]
        \unitlength=1mm
   \centering
   \begin{picture}(40,30)
      \put(1,2){\includegraphics[width=40\unitlength]{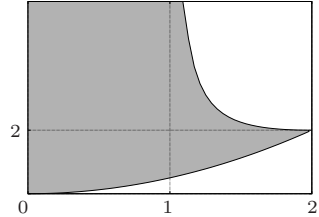}}
      \put(1,1){\scriptsize $0$}
      \put(20,1){\scriptsize $1$}
      \put(0,11){\scriptsize $2$}
      \put(39,1){\scriptsize $2$}
   \end{picture}
        \parbox[t]{0.80\textwidth}{\caption{\label{diagram-symmetric} All pairs $(w,A)$ satisfying the inequalities of Theorem~\ref{CentrSym}}}
\end{figure}

\begin{theorem} \label{CentrSym} 
Let $K \in \KK^2$ be lattice-free and centrally symmetric with
$\lwid:=\lwid(K)$ and $A:=\area(K)$.
Then 
\begin{align}
  0 & < \lwid \le 2, & & & & \label{SymLWidthUpper} \\[2mm]
  A & \le \infty & &\mbox{for} & 0 & < \lwid \le 1, \label{SymInf} \\ 
  A & \le \frac{\lwid^2}{2 (\lwid-1)} & &\mbox{for} & 1 & < \lwid
      \le 2, \label{SymAreaUpper} \\ 
  A & \ge \frac{1}{2} \lwid^2 & &\mbox{for} & 0&< \lwid \le
      2 \label{SymAreaLower}
\end{align}
(see Fig.~\ref{diagram-symmetric}). Furthermore, the following statements hold.
\begin{enumerate}[I.]
  \item The upper bound in \eqref{SymLWidthUpper} is attained if and
    only if, up to unimodular transformations, $$K = \conv \{ \pm
    (1,0), \pm (0,1) \} + \left(\frac{1}{2},\frac{1}{2}\right).$$
  \item \label{eq case for inf bound, sym} Equality in \eqref{SymInf}
    is attained if and only if $K$ is unbounded and contained in a
    split.
  \item \label{eq case for area bound} Equality in
    \eqref{SymAreaUpper} is attained if and only if, up to unimodular
    transformations,
    \[ K = \conv \left\{ \pm \left(\frac{\lwid}{2},0\right), \pm
       \left(0, \frac{\lwid}{2(\lwid-1)}\right) \right\} +
       \left(\frac{1}{2},\frac{1}{2}\right)\]
    (see Fig.~\ref{area-maximizer-sym-fig}).
  \item \label{eq case lower bound sym} Equality in
    \eqref{SymAreaLower} is attained if and only if, up to unimodular
    transformations, $K$ is a translate of
    \[ \frac{\lwid}{2}\conv \{ \pm (1,\alpha), \pm
       (0,1) \} \]
       for some $0 \le \alpha < 1$ satisfying $\max \{1+\alpha,2-\alpha\} \ge \lwid$ (see
       Fig.~\ref{area-minimizer-sym-fig}).
\end{enumerate}
\end{theorem}

%\begin{figure}[hbt]
%\begin{center}
%\unitlength=1mm
%\begin{picture}(70,70)
%\put(0,0){\includegraphics[width=40mm]{figures/sym-maximizer.eps}}
%\end{picture}
%\end{center}
%\caption{Maximizer of area for a given lattice width in the symmetric case}
%\end{figure}

%\begin{figure}[hbt]
%\begin{center}
%\unitlength=1mm
%\begin{picture}(70,70)
%\put(0,0){\includegraphics[width=40mm]{figures/sym-minimizer.eps}}
%\end{picture}
%\end{center}
%\caption{Minimizers of area for a given lattice width in the symmetric case}
%\end{figure}

A useful tool in convex geometry is the concept of covering minima
introduced by Kannan and Lov\'asz \cite{KannanLovasz88}.
We will therefore give an alternative formulation of our results in
terms of covering minima.
Let $K \in \KK^2$.
For $j = 1,2$ the $j$-th \term{covering minimum} is
\begin{equation*}
  \mu_j(K):= \inf \setcond{t \ge 0}{\text{each $(2-j)$-dimensional
  affine subspace of $\real^2$ intersects $t K + \integer^2$}}.
\end{equation*}

This definition implies $0 < \mu_1 \le \mu_2$.
We recall that $\mu_1(K) \wid(K) =1$.
Furthermore, for $t>0$, an appropriate translate of $t K$ is
lattice-free if and only if $t \le \mu_2(K)$.
Thus, Theorem~\ref{hurkens-thm} yields $\mu_2(K) \le
\left(1+\frac{2}{\sqrt{3}}\right) \mu_1(K)$.

For centrally symmetric bodies,
the lower bound in \eqref{SymAreaLower} was noticed by Makai
Jr.~\cite{Makai78}.
The upper bound in \eqref{SymLWidthUpper} is a consequence of a more
general result due to Kannan and Lov\'asz
\cite[Theorem\,2.13]{KannanLovasz88}.
Unfortunately, the proof of Theorem\,2.13 in \cite{KannanLovasz88} does not seem to
be correct\footnote{In the proof of Theorem\,2.13 in
\cite{KannanLovasz88} one claims that the covering minima of
centrally symmetric convex bodies satisfy the inequalities
$\mu_{k+1} \le 2 \mu_k$ (see \cite{KannanLovasz88} for the
explanation of the notations). In the proof of
Theorem\,2.13, p.~588, l.~11, it is inferred that
$2(\alpha + \beta + \lambda_1) \le 4(\alpha +
\beta)$. However, one line before it is just shown that
$\lambda_1 \le 2(\alpha + \beta)$ and therefore the correct
conclusion is $2(\alpha + \beta + \lambda_1) \le 6(\alpha +
\beta)$. Using the factor $6$ instead of $4$ in the rest of the proof
results in the weaker assertion $\mu_{k+1} \le 3\mu_k$. To the best of our
knowledge there is no revision or corrected version of this
proof which would yield the assertion $\mu_{k+1} \le 2 \mu_k$.}, but implies the weaker result $0 < \lwid \le 3$.
Therefore, we prove the subsequent theorem.

\begin{theorem} \label{cov min c-sym planar}
  Let $K \in \KK^2$ be centrally symmetric.
  Then $\mu_2(K) \le 2 \mu_1(K)$.
\end{theorem}

If $K$ is lattice-free, and thus $\mu_2(K) \ge 1$, then together with
the relation $\mu_1(K) \lwid(K) = 1$, Theorem \ref{cov min c-sym
  planar} implies $\lwid(K) \leq 2$.
The results stated in Theorems \ref{AreaLWidth} and \ref{CentrSym} can
also be expressed in terms of covering minima.
In Corollary \ref{Upper bounds for area in terms of mu's} the lower
bound for $A(K)$ in terms of $\mu_1(K)$ and $\mu_2(K)$ goes
back to \cite{MakaiToth74} (see also \cite{BetkeHenkWills93} for
further inequalities). 
Our results imply the sharp upper bounds for $A(K)$ in terms of
$\mu_1(K)$ and $\mu_2(K)$.

\begin{corollary} \label{Upper bounds for area in terms of mu's}
  Let $K \in \KK^2$ with $A:=A(K)$,  $\mu_1:=\mu_1(K)$ and
  $\mu_2:=\mu_2(K)$.
  Then
  \begin{align*}
	\mu_1 & \le \mu_2 \le \left( 1 + \frac{2}{\sqrt{3}} \right) \mu_1, \\
    A & \le \infty & &\mbox{for} &  \mu_1 & = \mu_2, \\
    A & \le \frac{1}{2 \mu_1 (\mu_2-\mu_1)} & &\mbox{for} & \mu_1 & <
        \mu_2 \le 2\mu_1, \\
    A & \le \frac{3 }{\mu_1 \left( 3 \mu_2 + \mu_1 - \sqrt{\mu_1^2 + 6
        \mu_1 \mu_2 - 3 \mu_2^2}\right) }&  &\mbox{for} & 2\mu_1 & <
        \mu_2 \le \left( 1 + \frac{2}{\sqrt{3}}\right) \mu_1, \\
    A & \ge \frac{3}{8\mu_1^2} & &\mbox{for} &  \mu_1 &  \le \mu_2 \le
        \left( 1 + \frac{2}{\sqrt{3}}\right) \mu_1.
  \end{align*}
  The upper bounds for $A$ are sharp, whereas the lower bound for $A$ is sharp only
  for $\mu_1 \le \mu_2 \leq 2 \mu_1$.
\end{corollary}

\begin{corollary} \label{mu.sym}
  Let $K \in \KK^2$ be centrally symmetric with $A:=A(K)$,
  $\mu_1:=\mu_1(K)$ and $\mu_2:=\mu_2(K)$.
  Then
  \begin{align*}
    \mu_1 & \le \mu_2 \le 2 \mu_1, & & & & \\[2mm]
    A & \le \infty & &\mbox{for} & \mu_1 & = \mu_2, \\
    A & \le \frac{1}{2 \mu_1 (\mu_2-\mu_1)} & &\mbox{for} & \mu_1 & <
        \mu_2 \le 2\mu_1, \\
    A & \ge \frac{1}{2\mu_1^2} & &\mbox{for} &  \mu_1 &  \le \mu_2 \le
        2 \mu_1.
  \end{align*}
  The above bounds are sharp.
\end{corollary}

For further inequalities between $\mu_1(K),\mu_2(K)$ and $\area(K)$
for the case $K \in \KK^2$ we refer to \cite{Schnell95}.

%%%%%%%%%%%%%%%%%%%%%%%%%%%%%%%%%%%%%%%%%%%%%%%%%%%%%%%%%%%%%%%%%%%%%%
%%%%%%%%%%%%%%%%%%%%%%%%%%%%%%%%%%%%%%%%%%%%%%%%%%%%%%%%%%%%%%%%%%%%%%

\section{Preliminaries} \label{prelims}

We consider the elements of $\real^2$ to be column vectors. The
transposition is denoted by $(\dotvar)^\top$ and the origin by $o$.
By $\aff$, $\conv$, $\bd$ and $\intr$ we denote the affine hull, the
convex hull, the boundary and the interior, respectively.
The maximum norm is denoted by $\|\dotvar\|_\infty$.
Remember that a set $K \in \KK^2$ is said to be \term{lattice-free} if
the interior of $K$ is disjoint with $\integer^2$.
A lattice-free set $K \in \KK^2$ is said to be
\term{maximal lattice-free} if $K$ is not properly contained in a
lattice-free set from $\KK^2$. The following fact is known.

\begin{proposition} \label{completing}
  Every lattice-free $K \in \KK^2$ is contained in a maximal
  lattice-free $H \in \KK^2$.
  %Furthermore, if $K$ is centrally symmetric, then $\Tilde{K}$
  %can be chosen to be centrally symmetric with the same center of
  %symmetry as $K$.
\end{proposition}

\newcommand{\UU}{\mathcal{U}}

\begin{proof}
For a closed convex set $U \subseteq \real^2$ and a point $x
\in \real^2$ we denote by $c(U,x)$ the topological closure
of the convex hull of $U \cup \{x\}$.

Let $(z_n)_{n=1}^{\infty}$ be the sequence of all elements
of $\rational^2$.
We define $U_0 := K$ and for every $n \in \natur$ we set
$U_n := c(U_{n-1},z_n)$ if $c(U_{n-1},z_n)$ is lattice-free,
and $U_n := U_{n-1}$ otherwise.
Let $H$ be the topological closure of
$\bigcup_{n=0}^{\infty}{U_n}$.
Since $U_{n-1} \subseteq U_n$ for every $n \in \natur$ it
holds $K \subseteq H$.
By construction, $H$ is a closed convex set with non-empty
interior, i.e., $H \in \KK^2$.
In addition, $H$ is lattice-free: assume $y$ is an integer
point in the interior of $H$.
Then there exists some $j \in \natur$ such that $y$ is in
the interior of $U_j$ and thus, $U_j$ is not lattice-free.
This contradicts the construction of $U_j$.
Let us show that $H$ is maximal lattice-free.
Assume the opposite and let $L \in \KK^2$ be lattice-free
such that $H \varsubsetneq L$.
Then $L \setminus H$ contains rational points which occur in
the sequence $(z_n)_{n=1}^{\infty}$.
Let $z_k \in L \setminus H$ be such a rational point. Since $U_{k-1} \subseteq H \subseteq L$ and $z_k \in L$, we have $c(U_{k-1},z_k) \subseteq L$. Thus, $c(U_{k-1},z_k)$ is lattice-free and, by definition of $U_k$, we have $U_k=c(U_{k-1},z_k)$. The latter implies $z_k \in H$, a contradiction to the choice of $z_k$.
\end{proof}

We point out that the above proof is not constructive. For a
constructive, but lengthy proof we refer to \cite{BCCZ09}.

A classification of planar maximal lattice-free convex sets was
noticed by Lov\'asz \cite{lovasz}. The refined classification which is
given below can be found in \cite{DeyWolsey}.

%\begin{proposition} \label{list-of-two-dim-max}
%  \thmcitation{\cite[Proposition\,3.3]{lovasz}}
%  Let $K \in \KK^2$ be maximal lattice-free.
%  Then $K$ is either a strip or a triangle or a quadrilateral and the
%  relative interior of each facet of $K$ contains at least one integer
%  point.
%\end{proposition}

%A classification of these sets can be found in \cite{DeyWolsey}.

\begin{figure}[ht]
  \unitlength=1mm
  \begin{center}
    \subfigure[\label{mlf.split} Split]{
        \begin{picture}(23,20)
                \put(2,0){\includegraphics[width=20\unitlength]{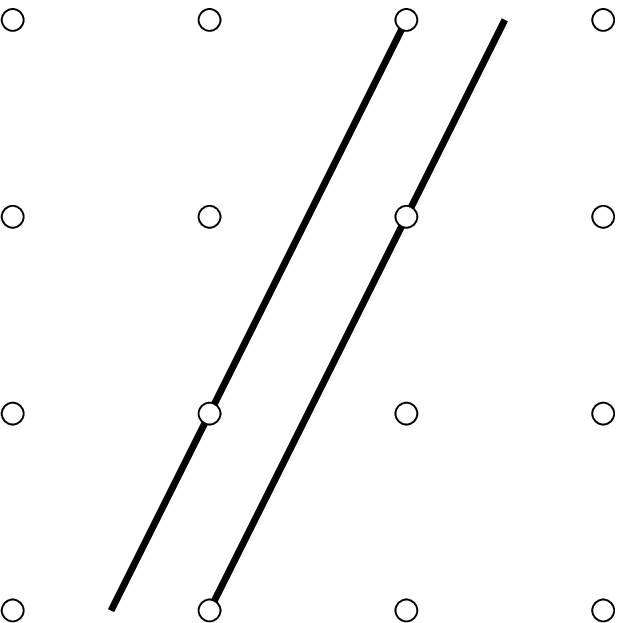}}
        \end{picture}
        } 
        \quad
    \subfigure[\label{mlf.tr1} Type 1 triangle]{
        \begin{picture}(23,20)
                \put(2,0){\includegraphics[width=20\unitlength]{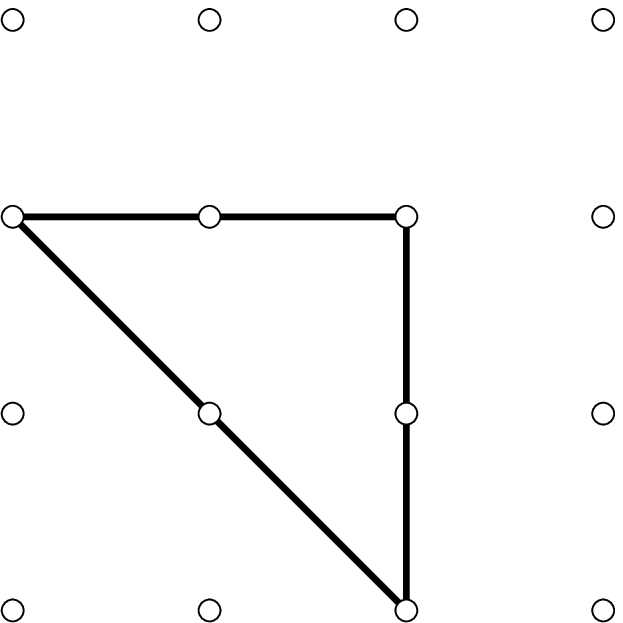}}
        \end{picture}
        } 
    \quad 
    \subfigure[\label{mlf.tr2} Type 2 triangle]{
        \begin{picture}(23,20)
                \put(2,0){\includegraphics[width=20\unitlength]{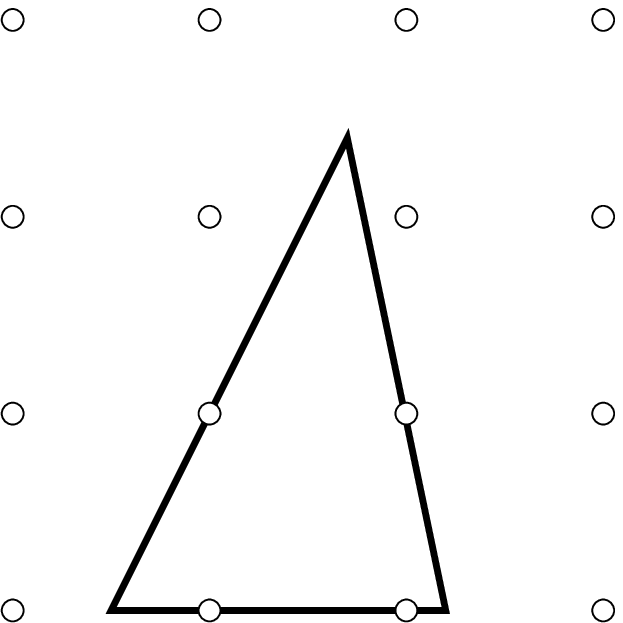}}
        \end{picture}
        }
    \quad
    \subfigure[\label{mlf.tr3} Type 3 triangle]{
        \begin{picture}(23,20)
                \put(2,0){\includegraphics[width=20\unitlength]{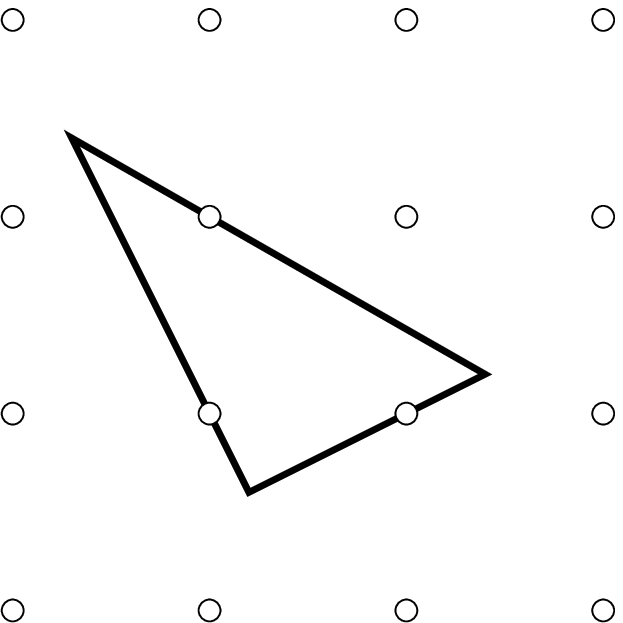}}
        \end{picture}
        } 
    \quad 
    \subfigure[\label{mlf.quad} Quadrilateral]{
        \begin{picture}(23,20)
                \put(2,0){\includegraphics[width=20\unitlength]{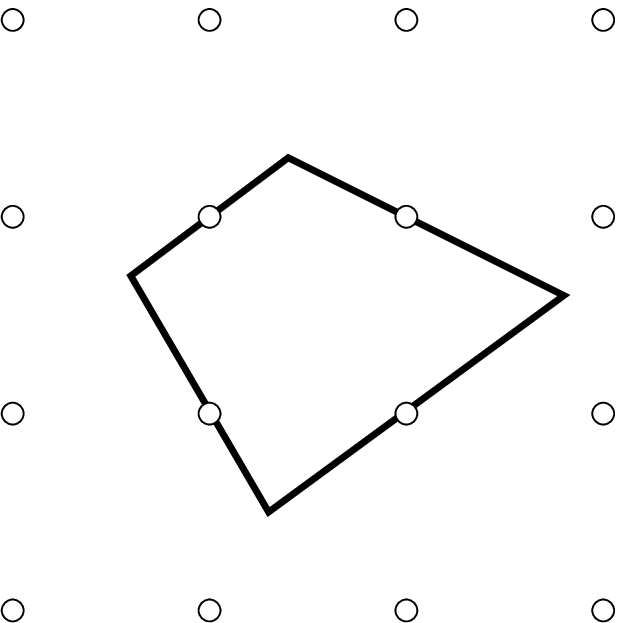}}
        \end{picture}
        }
  \end{center}
  \caption{\label{mlf.sets} All types of maximal lattice-free sets in $\KK^2$}
\end{figure}

\begin{proposition} \thmcitation{\cite{DeyWolsey} \cite{lovasz}}
  \label{list-of-two-dim-max}
Let $K \in \KK^2$ be maximal lattice-free.
Then $K$ is either a split or a triangle or a quadrilateral and the
relative interior of each facet of $K$ contains at least one integer
point.
In particular, $K$ is one of the following sets (see
Fig.~\ref{mlf.sets}).
\begin{enumerate} [I.]
  \item A split $\{(x_1,x_2) \in \real^2 : b \leq a_1 x_1 + a_2
    x_2 \le b + 1\}$ where $a_1$ and $a_2$ are coprime integers and
    $b$ is an integer.
  \item A triangle with at least one integer point in the relative
    interior of each of its edges, which in turn is either
    \begin{enumerate} [(a)]
      \item a type 1 triangle, i.e., a triangle with integer vertices
        and exactly one integer point in the relative interior of
        each edge, or
      \item a type 2 triangle, i.e., a triangle with at least one
        fractional vertex $v$, exactly one integer point in the
        relative interior of the two edges incident to $v$ and at
        least two integer points on the third edge, or
      \item a type 3 triangle, i.e., a triangle with exactly three
        integer points on the boundary, one in the relative interior
        of each edge.
    \end{enumerate}
  \item A quadrilateral containing exactly one integer point in the
    relative interior of each of its edges.
\end{enumerate}
\end{proposition}

Theorems \ref{AreaLWidth} and \ref{CentrSym} will be proved for
maximal lattice-free convex sets first.
Once this is established we show that this implies their validity for
lattice-free convex sets which are not maximal as well.

For $K \in \KK^2$, the \term{support function} of $K$ is defined by
$h(K,u):=\sup \setcond{u^\top x}{x \in K}$.
Remember that the \term{width function} of $K$ is
$\wid(K,u) = h(K,u)+h(K,-u)$, where $u \in \real^2$, and the
\term{lattice width} of $K$ is the value $\lwid(K) = \min
\setcond{w(K,u)}{u \in \integer^2 \setminus \{o\}}$.
Note that the lattice width is invariant with respect to unimodular
transformations. Bounded elements of $\KK^2$ are referred to as \term{convex bodies}.
If $K  \in \KK^2$ is a convex body, then $DK:=\setcond{x-y}{x, y \in
  K}$ is said to be the \term{difference body} of $K$.
It is well known that $h(DK,u) = w(K,u)$ for every $u \in \real^2$.
If $K$ is a convex body symmetric in the origin, then the
\term{Minkowski functional} $\| \dotvar \|_K$ of $K$ is defined by
$\|u\|_K := \min \setcond{\lambda \ge 0}{u \in \lambda K}$, where $u
\in \real^2$.
For a convex body $K$ containing the origin in the interior the set
$K^\ast := \setcond{u \in \real^2}{h(K,u) \le 1}$ is a convex body as
well and is referred to as the \term{polar body} of $K$.
One has $h(K,u) = \|u\|_{K^\ast}$ for every $u \in \real^2$.
The following formula, in slightly different terms, can be found in
\cite[Lemma~2.3]{KannanLovasz88}:
\begin{equation}
  \lwid(K) = \sup \setcond{ \alpha>0}{ (\alpha  (DK)^\ast ) \cap
  \integer^2 = \{o\}}. \label{lattice-width-over-lambda}
\end{equation}

A collection of convex sets $\mathcal{X}$ from $\KK^2$ is said to
\term{tile} $\real^2$ if the interiors of the elements of
$\mathcal{X}$ are pairwise disjoint and the union of $\mathcal{X}$
yields $\real^2$ (see also \cite[Section~4.1]{Schulte93} for
information on lattice tilings).

The two-dimensional version of \term{Minkowski's first fundamental
theorem} (see \cite[Theorem\,22.1]{GruberBook07}) states that if $K \in
\KK^2$ is symmetric in the origin and $\intr K \cap \integer^2 =
\{o\},$ then $\area(K) \le 4$.
Furthermore, if $\area(K)=4$, then the sets $\frac{1}{2} K + z$ with
$z \in \integer^2$ tile $\real^2$.
\term{Mahler's inequality} (see \cite{Mahler39}) states that $\area(K)
\area(K^\ast) \ge 8$ for convex bodies $K \in\KK^2$ symmetric in the
origin, with equality if and only if $K$ is a parallelogram. 
The following proposition is easy to show.

\begin{proposition} \label{tiling-by-parallelograms}
  Let $P$ be a parallelogram symmetric in the origin and such that its
  translates $P+z$ with $z \in \integer^2$ tile $\real^2$.
  Then, up to unimodular transformations, $P= \frac{1}{2} \conv \{\pm(-\alpha-1,1),\pm(-\alpha+1,1)\}$ for some $0 \le \alpha < 1$.
\end{proposition}

Proposition~\ref{tiling-by-parallelograms} is a special case
of a result of Haj\'os \cite[\S\S1,2]{Hajos41}, see also
\cite[p.\,174]{GruLek87}.

%%%%%%%%%%%%%%%%%%%%%%%%%%%%%%%%%%%%%%%%%%%%%%%%%%%%%%%%%%%%%%%%%%%%%%
%%%%%%%%%%%%%%%%%%%%%%%%%%%%%%%%%%%%%%%%%%%%%%%%%%%%%%%%%%%%%%%%%%%%%%

\section{Auxiliary results on triangles} \label{aux:tri}

An essential part of the proofs of our main results is concerned with
analytical representations of the lattice width and the area of
various types of maximal lattice-free polygons.
In particular, we shall need formulas for the lattice width of
triangles.

We start the section by presenting well-known facts about barycentric
coordinates.
\begin{lemma} \label{BaryLem}
Let $q_0,q_1,q_2$ be affinely independent points in $\real^2$ and let
$p, p_0, p_1, p_2$ be points in $\real^2$ which are represented in the
form
\begin{align*}
  p &=\sum_{j=0}^2 x_j q_j & &\mbox{and} & p_i & = \sum_{j=0}^2
  x_{i,j} q_j,
\end{align*}
where
\begin{align*}
  1 & = \sum_{j=0}^2 x_j & &\mbox{and} & 1 = \sum_{j=0}^2 x_{i,j},
\end{align*}
and $x_j, x_{i,j} \in \real$ for $i,j = 0,1,2$. We define $l_j:= \aff
(\{q_0,q_1,q_2\} \setminus \{q_j\})$ for
$j=0,1,2$.
Then the following statements hold.
\begin{enumerate}[I.]
  \item \label{BarySameSide} The points $p$ and $q_j$ lie in the same open halfplane
    determined by $l_j$ if and only if $x_j>0$.
  \item The value $|x_j|$ is the ratio of the distance from $p$ to
    $l_j$ and the distance from $q_j$ to $l_j$.
  \item \label{BaryArea} The areas of $Q:=\conv \{q_0,q_1,q_2\}$ and
    $P:=\conv \{p_0,p_1,p_2\}$ are related by
    \begin{equation*}
      \area(P) = \big| \det (x_{i,j})_{i,j=0,\dots,2} \big| \area(Q).
    \end{equation*}
\end{enumerate}
\end{lemma}
The values $x_0,x_1,x_2$ associated to a point $p$ in
Lemma~\ref{BaryLem} are said to be the \term{barycentric coordinates}
of $p$ with respect to the triangle $\conv \{q_0,q_1,q_2\}$, see
\cite[Sect.~13.7]{Coxeter89}.

\newcommand{\MatrixOfTriangle}[1]{
        \begin{pmatrix}
                {#1}_0^\top & 1 \\
                {#1}_1^\top & 1 \\
                {#1}_2^\top & 1 
        \end{pmatrix}   
}
\newcommand{\qMx}{\MatrixOfTriangle{q}}
\newcommand{\pMx}{\MatrixOfTriangle{p}}
\newcommand{\dMx}{
        \begin{pmatrix*}[r]
                -1 & 1 & 0 \\
                0 & -1 & 1 \\
                1 & 0 & -1 
        \end{pmatrix*}
}

\begin{lemma} \label{tri-wid-analyticially}
Let $P:=\conv \{p_0,p_1,p_2\}$ such that $p_0,p_1,p_2 \in \integer^2$
are the only integer points in $P$.
Let $Q:=\conv \{q_0,q_1,q_2\}$ be a triangle whose vertices are given
by the barycentric coordinates with respect to $P$, that is, by a $3
\times 3$-matrix $B$ such that
\begin{equation*}
        \qMx =  B \pMx
\end{equation*}
Then
\begin{equation} \label{tri-wid-formula-1}
        \lwid(Q) = \min \setcond{\|D B z\|_\infty}{z \in \integer^3 \ \mbox{and the coordinates of $z$ are not all equal}},
\end{equation}
where 
\begin{equation*}
        D:=\dMx.
\end{equation*}

Furthermore, if $p_i:=(1-x_i) q_{i+1} + x_i q_{i+2}$ with $i=0,1,2$
and $0<x_i < 1$ (that is, $Q$ is circumscribed about $P$), then
% \begin{equation} \label{tri-wid-formula-2}
%   \lwid(Q) = \frac{1}{x_0x_1x_2 +
%     (1-x_0)(1-x_1)(1-x_2)} \min\limits_{\overtwocond{y \in
%     \integer^3 \setminus\{o\}}{y_0+y_1+y_2=0}}  \max_{i=0,1,2}\limits
%     |x_i y_i + (1-x_{i+1}) y_{i+1}| 
% \end{equation}
 \begin{equation} \label{tri-wid-formula-2}
   \lwid(Q) = \frac{\min \setcond{ \max_{i=0,1,2}\limits
     |x_i y_i + (1-x_{i+1}) y_{i+1}| }{y \in
     \integer^3 \setminus\{o\}, \ y_0+y_1+y_2=0}}{x_0x_1x_2 +
     (1-x_0)(1-x_1)(1-x_2)}
 \end{equation}
and 
\begin{equation} \label{area-formula}
        \area(Q) = \frac{1}{2 (x_0x_1x_2 + (1-x_0)(1-x_1)(1-x_2)) }.
\end{equation}
\end{lemma}

\begin{proof}
  For $u \in \integer^2$ we have
  \begin{align*}
    w(Q,u) & = \max \setcond{|q_i^\top u - q_j^\top u|}{0\le i < j \le
      2} = \left\| D
      \begin{pmatrix}q_0^\top \\ q_1^\top \\ q_2^\top\end{pmatrix} u
      \right\|_\infty = \left\| D \qMx
      \begin{pmatrix} u \\ k \end{pmatrix} \right\|_\infty \\
      & = \left\| D B \pMx \begin{pmatrix} u \\ k \end{pmatrix}
      \right\|_\infty = \left\| D B z \right\|_\infty,
  \end{align*}
where $k \in \integer$ is arbitrary and 
\begin{equation*}
  z:=\pMx \begin{pmatrix} u \\ k \end{pmatrix} = \begin{pmatrix}
  p_0^\top u + k \\ p_1^\top u + k \\ p_2^\top u + k\end{pmatrix}.
\end{equation*}
Clearly, $z \in \integer^3$.
Since the vector $z$ is the product of a unimodular matrix and an
integer vector, it follows that $u=o$ if and only if the coordinates
of $z$ are all equal. This shows \eqref{tri-wid-formula-1}.
                
Let us show \eqref{tri-wid-formula-2}. We have
\begin{align*}
  & \pMx =
  X \qMx, &
  & \mbox{where} &
  X := B^{-1} = \begin{pmatrix}
                  0 & 1-x_0 & x_0 \\ 
                  x_1 & 0 & 1-x_1 \\
                  1-x_2 & x_2 & 0
                \end{pmatrix}.
\end{align*}
Since $X$ is the matrix of barycentric coordinates of the vertices of
$P$ with respect to $Q$, $B$ is the matrix of barycentric coordinates
of the vertices of $Q$ with respect to $P$.
Direct computations yield
\begin{align*}
  \det X & = x_0 x_1 x_2 + (1-x_0) (1-x_1) (1-x_2)>0, \\
  B & = \frac{1}{\det X}
  \begin{pmatrix} -(1-x_1)x_2 & x_0 x_2 & (1-x_0)(1-x_1) \\
    (1-x_1)(1-x_2) & - (1-x_2) x_0 & x_0 x_1 \\
    x_1 x_2 & (1-x_0)(1-x_2) & -(1-x_0) x_1
  \end{pmatrix}, \\
  D B & = \frac{1}{\det X} 
  \begin{pmatrix}
        1-x_1 & -x_0 & x_0+x_1-1 \\
        x_1+x_2 - 1 & 1-x_2 & -x_1 \\
        -x_2 & x_0+x_2 -1 & 1-x_0
  \end{pmatrix}.
\end{align*}
We employ the latter matrix relation and obtain for every $z \in
\integer^3$ the following:
\begin{align*}
        D B z  = \frac{1}{\det X}
        \begin{pmatrix} 
                z_2-z_1 & z_2-z_0 & 0 & z_0-z_2 \\
                0 & z_0-z_2 & z_0-z_1 & z_1-z_0 \\
                z_1-z_2 & 0 & z_1-z_0 & z_2-z_1
        \end{pmatrix}
        \begin{pmatrix}
                x_0 \\
                x_1 \\
                x_2 \\
                1
        \end{pmatrix}
\end{align*}
By the change of variables
\begin{equation} \label{y:z:change}
        y := \begin{pmatrix}z_2 - z_1 \\ z_0-z_2 \\z_1 - z_0\end{pmatrix}
\end{equation}
the latter amounts to 
\begin{align*}
  D B z  & = \frac{1}{\det X}
  \begin{pmatrix*}[r] 
    y_0 & -y_1 & 0 & y_1 \\
    0 & y_1 & -y_2 & y_2 \\
    -y_0 & 0 & y_2 & y_0
  \end{pmatrix*}
  \begin{pmatrix}
    x_0 \\
    x_1 \\
    x_2 \\
    1
  \end{pmatrix}
  \\
  & = \frac{1}{\det X} \begin{pmatrix}x_0 y_0 + (1-x_1) y_1 \\ x_1 y_1
  + (1-x_2) y_2 \\ x_2 y_2 + (1-x_0) y_0 \end{pmatrix}. 
\end{align*}
Clearly, $y_0+y_1+y_2=0$ and, if the coordinates of $z$ are not all
equal, we have $y \ne o$.
Conversely, for an arbitrary $y \in \integer^3 \setminus \{o\}$ with
$y_0+y_1+y_2=0$, we can easily find an appropriate $z \in \integer^3$
satisfying \eqref{y:z:change}.
Thus, employing \eqref{tri-wid-formula-1} and the previous derivation
we arrive at
\begin{equation*}
  \lwid(Q) = \frac{1}{\det X} \min \setcond{ \max
    \setcond{|x_i y_i + (1-x_{i+1}) y_{i+1}|}{i=0,1,2} }{y \in
    \integer^3 \setminus\{o\}, \ y_0+y_1+y_2=0}.
\end{equation*}

Equation \eqref{area-formula} is a straightforward consequence of Lemma~\ref{BaryLem}.\ref{BaryArea}.
\end{proof}

\section{Proofs for arbitrary bodies} \label{general.case}

In order to prove Theorem \ref{AreaLWidth} it is convenient to show
the statement only for maximal lattice-free convex
sets, instead of lattice-free ones.
This is possible since every lattice-free convex set is contained in a
maximal lattice-free convex set, see Proposition \ref{completing}.

Lemma~\ref{FreeCond} characterizes, in analytic terms, maximal
lattice-free triangles of type 3 and their lattice width.

\begin{figure}[hbt]
\begin{center}
\unitlength=0.6mm
\begin{picture}(70,70)
        \put(5,5){\includegraphics[width=60\unitlength]{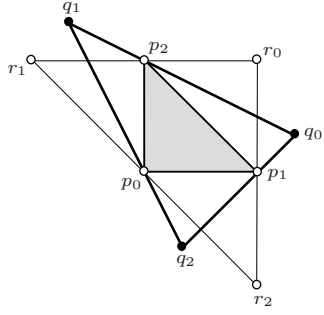}}
        \put(26,28){\scriptsize $p_0$}
        \put(58,30){\scriptsize $p_1$}
        \put(32,58){\scriptsize $p_2$}
        \put(66,39){\scriptsize $q_0$}
        \put(13,67){\scriptsize $q_1$}
        \put(38,11){\scriptsize $q_2$}
        \put(57,57){\scriptsize $r_0$}
        \put(1,53){\scriptsize $r_1$}
        \put(55,2){\scriptsize $r_2$}
%       \graphpaper[5](0,0)(70,70)
\end{picture}
\end{center}
\caption{\label{Fig1} Points $p_i, q_i, r_i$, $i\in \{0,1,2\}$,
as in the proof of Lemma~\ref{FreeCond}}
\end{figure}

\begin{lemma} \label{FreeCond}
Let $P:=\conv \{p_0,p_1,p_2\}$ such that $p_0,p_1,p_2 \in \integer^2$
are the only integer points in $P$.
Let $Q:=\conv \{q_0,q_1,q_2\}$ be a triangle circumscribed about $P$
so that $p_i:=(1-x_i) q_{i+1} + x_i q_{i+2}$ for $i=0,1,2$ and
$0<x_i<1$ (see also Fig.\,\ref{Fig1}).
Then the following statements hold.
\begin{enumerate}[I.]
  \item \label{LatFreeCond} $Q$ is a maximal lattice-free triangle of
    type 3 if and only if
    \begin{enumerate}[(a)]
      \item $x_i+x_j > 1$ for all $0 \le i < j \le 2$ or
      \item $x_i+x_j < 1$ for all $0 \le i < j \le 2$.
    \end{enumerate}
  \item \label{LWidthFormula} If $(a)$ holds, then the lattice width
    of $Q$ is given by
    \begin{equation*}
      \lwid(Q) = \frac{\min \{x_0,x_1,x_2\}}{x_0 x_1 x_2 + (1-x_0)
      (1-x_1) (1-x_2)}.
    \end{equation*}
  \item \label{LWidthBound} If $(a)$ holds, then $\lwid(Q) \le
    1+\frac{2}{\sqrt{3}}$ with equality if and only if $x_0=x_1=x_2 =
    \frac{1}{\sqrt{3}}$.
\end{enumerate}
%Then $Q$ is a maximal lattice-free triangle of type 3 if and only if
%    \begin{enumerate}[(a)]
%      \item $x_i+x_j > 1$ for all $0 \le i < j \le 2$ or
%      \item $x_i+x_j < 1$ for all $0 \le i < j \le 2$.
%    \end{enumerate}
%
%Furthermore, if (a) is fulfilled, one has $\lwid(Q) \le
%    1+\frac{2}{\sqrt{3}}$ with equality if and only if $x_0=x_1=x_2 =
%    \frac{1}{\sqrt{3}}$.
\end{lemma}

\begin{proof}
Assume that $Q$ is a lattice-free triangle of type 3.
By $H_i$ we denote the closed halfplane with $q_{i+1},q_{i+2} \in \bd
H_i$ and $q_i \in H_i$.
We also introduce the points $r_i:=-p_i + p_{i+1} + p_{i+2} \in
\integer^2$.
By construction, $p_i$ is the midpoint of $\conv \{r_{i+1},r_{i+2}\}$.
Because of the latter property, and since $p_i \in \bd H_i$, we have
$r_{i+1} \in H_i$ or $r_{i+2} \in H_i.$
For $i \in \{0,1,2\}$ by $\tau(i)$ we denote the set of all $k \in
\{0,1,2\}$ such that $k \ne i$ and $r_k \in H_i.$
By the above observations $\tau(i) \ne \emptyset$ for every $i$.
If for some $0 \le i < j \le 2$ one has $\tau(i) \cap \tau(j) \ne
\emptyset$ we choose $k \in \tau(i) \cap \tau(j)$.
Then $r_k \in H_0 \cap H_1 \cap H_2 = Q$, and by this the cell $\conv
\{p_0,p_1,p_2,r_k\}$ of $\integer^2$ is a subset of $Q.$
The latter means that $Q$ is not of type 3, a contradiction.
Thus, $\tau(i) \cap \tau(j) = \emptyset$ for $0 \le i < j \le 2$.
Taking into account that $i\not\in \tau(i)$ for $i \in \{0,1,2\}$, we
see that $\tau(i)$ is a singleton for every $i$ and is in fact one of
the two possible cyclic shifts on $\{0,1,2\}$.
In other words, either $\tau(i) = i+1 \modulo{3}$ for every $i$ or
$\tau(i)=i+2 \modulo{3}$ for every $i$.
If $\tau(i)=i+1 \modulo{3}$ for every $i$, then $r_{i+2} \not\in H_i$
for every $i.$
This means, that the $i$-th barycentric coordinate of $r_{i+2}$ with
respect to $Q$ is strictly negative.
This barycentric coordinate is $x_{i+1}+x_{i+2}-1$ since 
\begin{align*}
        r_{i+2} & =-p_{i+2}+p_i+p_{i+1} \\ & = -((1-x_{i+2})q_i+x_{i+2}q_{i+1}) + ((1-x_i)q_{i+1}+x_i q_{i+2}) + ((1-x_{i+1}) q_{i+2} + x_{i+1} q_i) \\
                & = (x_{i+1}+x_{i+2}-1) q_i + (1-x_i+x_{i+2}) q_{i+1} + (1-x_{i+1}+x_i) q_{i+2}.
\end{align*}
Thus, we obtain (b).
If $\tau(i)=i+2 \modulo{3}$ for every $i$, arguing in the same way we
obtain (a).

For proving the converse, we assume that (a) or (b) is fulfilled and
show that $Q$ is a lattice-free triangle of type 3.
Consider an arbitrary $p \in \integer^2.$
We can represent $p$ by $p=z_0 p_0 + z_1 p_1 + z_2 p_2$ where $z_i \in
\integer$ and $z_0+z_1+z_2=1.$
From symmetry reasons, we may assume that $z_0 \le z_1 \le z_2$.
Under these assumptions, we have $z_2 \ge 1$ and $z_0 \le 0$.
We evaluate the barycentric coordinates of $p$ with respect to $Q$ as
follows:
\begin{align*}
        p & = z_0 p_0 + z_1 p_1 + z_2 p_2 \\
        & = z_0 ((1-x_0)q_1 + x_0 q_2) +z_1 ((1-x_1)q_2 + x_1 q_0)
        + z_2 ((1-x_2)q_0 + x_2 q_1) \\
        & = (z_1 x_1 + z_2 (1-x_2)) q_0 + (z_0(1-x_0)+z_2 x_2) q_1
        + (z_0 x_0 + z_1 (1-x_1)) q_2. 
\end{align*}
If $z_1 \le 0$, then the barycentric coordinate $z_0 x_0 + z_1
(1-x_1)$ is non-positive and by this, in view of Lemma~\ref{BaryLem}.\ref{BarySameSide}, $p \not\in \intr Q.$
Assume that $z_1 \ge 1$.
If (a) is fulfilled, then the barycentric coordinate $z_0 x_0 + z_1
(1-x_1)$ is estimated as follows:
\begin{align*}
  z_0 x_0 + z_1 (1-x_1) < z_0 x_0 + z_1 x_0 = (z_0+z_1) x_0 =
  (1 - z_2) x_0 \le 0.
\end{align*}
Consequently, $p \not\in \intr Q$.
If (b) is fulfilled, the barycentric coordinate $z_0 (1-x_0)+z_2 x_2$
can be estimated analogously:
\begin{align*}
  z_0 (1-x_0)+z_2 x_2 < z_0 (1-x_0) + z_2 (1-x_0) = (z_0+z_2) (1-x_0)
  = (1-z_1) (1-x_0) \le 0.
\end{align*}
%If (b) is fulfilled, it follows analogously:
%\begin{align*}
%  z_0 x_0 + z_1 (1-x_1) < z_0 (1-x_1) + z_1 (1-x_1) =
%  (z_0+z_1) (1-x_1) = (1-z_2) (1-x_1) \le 0.
%\end{align*}
Thus, also in this case $p \not\in \intr Q$.
This shows the first part of the lemma.

We now show the second part. In view of Lemma~\ref{tri-wid-analyticially} it suffices to show
$$g(x):=\min \setcond{ \max_{i=0,1,2}\limits |x_i y_i + (1-x_{i+1})
  y_{i+1}| }{y \in \integer^3 \setminus\{o\}, \ y_0+y_1+y_2=0}
  = \min \{x_0,x_1,x_2\}$$
under the assumption that $x_i+x_j > 1$ for all $0 \le i < j \le 2$.
Taking all six choices of $y \in \{-1,0,1\}^3$ with $y_0+y_1+y_2=0$
we easily verify that 
$$g(x) \le \min_{i=0,1,2} \max \{x_{i+1},1-x_i\}
  = \min \{x_0,x_1,x_2\},$$
where the last inequality is due to assumption (a).
It remains to show the converse inequality.
Consider $y \in \integer^3 \setminus \{o\}$ with $y_0+y_1+y_2=0$.
Possibly, interchanging $y$ by $-y$, we assure the existence of $j$
such that $y_j \ge 0$ and $y_{j+1} \ge 0.$
If $y_j \ge 1$ and $y_{j+1} \ge 1$, then 
%\begin{align*}
%  \max_{i=0,1,2} |x_i y_i + (1-x_{i+1}) y_{i+1}| & \ge
%  x_j+ 1-x_{j+1} \ge \max \{x_j,1-x_{j+1}\} \\ 
%  & \ge \min_{i=0,1,2} \max \{x_i,1-x_{i+1}\} \\
%  & = \min \{x_0,x_1,x_2\}.
%\end{align*}
\begin{align*}
  \max_{i=0,1,2} |x_i y_i + (1-x_{i+1}) y_{i+1}| & \ge
  x_j+ 1-x_{j+1} \ge x_j \ge \min \{x_0,x_1,x_2\}.
\end{align*}
Otherwise, one of the $y_i$'s is equal to zero and the remaining ones
are equal to $k$ and $-k$ for some $k \in \natur$.
The latter means we can replace $y$ by $\frac{1}{k} y,$ which would
decrease $\max \setcond{|x_i y_i + (1-x_{i+1}) y_{i+1}|}{i=0,1,2}$.
We thus arrive at the case $y \in \{-1,0,1\}^3$ which has already been
considered above.
%This shows the converse inequality and proves
%\eqref{tri-wid-formula-2}.

Let us prove the third part of the lemma.
Without loss of generality let $x_0 \le x_1 \le x_2$.

\emph{Case 1:} $x_0 \le \frac{1}{2}$.
Then $x_1 > \frac{1}{2}$ and $x_2 > \frac{1}{2}$ and we have
\begin{align}
  \frac{1}{\lwid(Q)} & = 
  \frac{x_0 x_1 x_2 + (1-x_0) (1-x_1) (1-x_2)}{x_0} \nonumber \\
  & \ge \frac{x_0 x_1 x_2 + x_0 (1-x_1) (1-x_2)}{x_0} \nonumber \\
  & = \frac{1}{2} (2 x_1 - 1) (2x_2 - 1) + \frac{1}{2}
  > \frac{1}{2}, \label{x0<=1 implies w<=2}
\end{align}
which implies that $\lwid(Q) < 2$.

\emph{Case 2:} $x_0> \frac{1}{2}$.
We use the notations $\sigma_1(x) :=x_0+x_1+x_2$ and $\sigma_2(x):=x_0
x_1+x_0 x_2 + x_1 x_2$.
Clearly, $2 \sigma_2(x) - \sigma_1(x) = (x_0+x_1-1) x_2+ (x_0+x_2-1)
x_1 + (x_1+x_2-1) x_0 \ge (2 (x_0+x_1+x_2)-3) x_0 = (2 \sigma_1(x) -
3)x_0$ and by this 
\begin{align}
  x_0 x_1 x_2 + (1-x_0) (1-x_1) (1-x_2) & = 1-\sigma_1(x) + \sigma_2(x) \ge 1- \frac{1}{2}
  \sigma_1(x) + \left(\sigma_1(x) - \frac{3}{2}\right)x_0 \nonumber \\
  & = 1-\frac{3}{2} x_0  + \left(x_0 - \frac{1}{2}\right) \sigma_1(x)
  \label{detA lower} \\
  & \ge 1- \frac{3}{2} x_0 + \left(x_0-\frac{1}{2}\right) 3 x_0 = 1 -
  3 x_0 + 3 x_0^2. \nonumber
\end{align}
Consequently, applying elementary calculus, we get
\begin{equation*}
  \frac{1}{\lwid(Q)} \ge 3 x_0 -3 +
  \frac{1}{x_0} \ge 2 \sqrt{3} - 3,
\end{equation*}
which implies $\lwid(Q) \le 1+\frac{2}{\sqrt{3}}$ and shows that
the equality $\wid(Q)=1+ \frac{2}{\sqrt{3}}$ is attained if and
only if $x_0=x_1=x_2=\frac{1}{\sqrt{3}}$.
\end{proof}

The following two lemmas prepare the proof of Theorem
\ref{AreaLWidth}.
Parts of the proof of Lemma \ref{bounds:for:quad} are borrowed
from \cite{Hurkens90}; nevertheless we need these parts for subsequent
arguments.

\begin{lemma} \label{bounds:for:quad}
  Let $K \in \KK^2$ be maximal lattice-free with $[0,1]^2 \subseteq K$
  and let $\lwid:=\lwid(K)$ and $A:=\area(K)$.
  Then $\lwid \le 2$ and either it holds $\lwid = 1$ and $A = \infty$
  (i.e.~$K$ is a split) or $\lwid > 1$ and $A \le \frac{\lwid^2}{2
  (\lwid -1 )}$ with equality $A= \frac{\lwid^2}{2 (\lwid -1 )}$
  characterized by Part~\ref{AreaEasyEq} of Theorem~\ref{AreaLWidth}.
\end{lemma}

\begin{figure}[hbt]
\begin{center}
\unitlength=0.7mm
\begin{picture}(50,60)
\put(0,0){\includegraphics[width=50\unitlength]{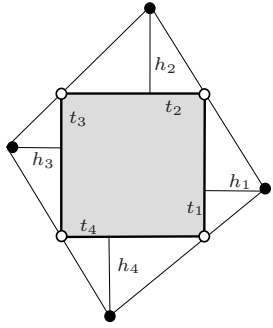}}
\put(42,26){\scriptsize$h_1$}
\put(28,48){\scriptsize$h_2$}
\put(5,30){\scriptsize$h_3$}
\put(21,10){\scriptsize$h_4$}
\put(34,21){\scriptsize$t_1$}
\put(30,40){\scriptsize$t_2$}
\put(12,38){\scriptsize$t_3$}
\put(14,17.5){\scriptsize$t_4$}
\end{picture}

\caption{\label{QuadFig} A maximal lattice-free
  quadrilateral in the proof of Theorem~\ref{AreaLWidth}}
\end{center}
\end{figure}

\begin{proof}
If $\lwid = 1$ there is nothing to show.
Thus, assume $\lwid > 1$ and therefore $K$ is a triangle of type 1,
type 2 or a quadrilateral.
We only consider the case that $K$ is a quadrilateral.
The case where $K$ is a triangle can be viewed as a degenerate version
of a quadrilateral where one vertex becomes a convex combination of
its two neighbor vertices.
By $a_1,a_2,a_3, a_4$ we denote the consecutive vertices of $[0,1]^2$.
Let $q_1,q_2, q_3, q_4$ be consecutive vertices of $K$ such that the
point $q_i'$ of $[0,1]^2$ closest to $q_i$ lies in $\conv
\{a_i,a_{i+1}\}$.
The distance from $q_i$ to $q_i'$ will be denoted by $h_i$ and the
distance from $a_i$ to $q_i'$ by $t_i$, see also Fig.~\ref{QuadFig}.
Taking into account the relations $h_i h_{i-1} = t_i (1-t_{i-1})$ it
can be verified that
\begin{align}
    &~1 - (h_1+h_3)(h_2+h_4) \nonumber \\
  = &~(1-t_1-t_3)(1-t_2-t_4) \label{hbounds} \\
  = &~\frac{(h_1h_2h_3h_4-t_1t_2t_3t_4)^2}{t_1 t_2 t_3
        t_4} \ge 0, \nonumber
\end{align}
for details see \cite[p.\,124]{Hurkens90}.
Without loss of generality we assume that $h_1+h_3 \le h_2 + h_4$.
Relations \eqref{hbounds} yield $h_1+h_3 \le 1$.
Thus, the width of $K$ with respect to the vector $u = (1,0)$ is
$h_1+h_3 + 1 \le 2$.
For all vectors $u \in \integer^2 \setminus \{o\}$ which are not in
$\{\pm (1,0), \pm (0,1)\}$ we easily get $\wid(K,u) \ge
\wid([0,1]^2,u) \ge 2$.
Hence $\lwid = h_1+h_3+1 \le 2$. Furthermore,
\begin{align*}
  A & = 1 + \frac{1}{2} (h_1+h_2+h_3+h_4) = 1 + \frac{1}{2} (\lwid -1
  + h_2 +h_4) \\
  & \le 1 + \frac{1}{2} \left(\lwid -1  + \frac{1}{h_1+h_3} \right) =
  1 + \frac{1}{2} \left(\lwid - 1 + \frac{1}{\lwid -1}\right) =
  \frac{\lwid^2}{2(\lwid-1)}.
\end{align*}
If the equality $A= \frac{\lwid^2}{2 (\lwid -1 )}$ is attained, then
$h_2+h_4+1=\frac{\lwid}{\lwid-1}$ and $(1-t_1-t_3)(1-t_2-t_4)=0,$
which implies that $1-t_1-t_3=0$ or $1-t_2-t_4=0$.
Taking into account the geometric meaning of $t_i$'s, we see (using
\eqref{hbounds}) that the equalities $1-t_1-t_3=0$ and $1-t_2-t_4=0$
imply one another so that one has $1-t_1-t_3=1-t_2-t_4=0$.
The above relations yield the characterization of the equality case
given in Part~\ref{AreaEasyEq} of Theorem~\ref{AreaLWidth}.
\end{proof}

\begin{lemma} \label{bounds:for:tri}
  Let $K \in \KK^2$ be a maximal lattice-free triangle with
  $\lwid:=\lwid(K)$ and $A:=\area(K)$.
  Then $\lwid>1$ and \eqref{AreaEasy} resp.~\eqref{AreaHard} holds
  true.
  The equality case in both inequalities is characterized by
  Part~\ref{AreaEasyEq} resp.~Part \ref{AreaHardEq} of
  Theorem~\ref{AreaLWidth}.
\end{lemma}

\begin{proof}
If $K$ is a triangle of type 1 or type 2, then clearly $\lwid > 1$.
Consider a triangle $K$ of type 3 and a set $P = \conv
\{p_0,p_1,p_2\}$, as in Lemma \ref{FreeCond}, such that the relative
interior of each edge of $K$ contains a point from $\{p_0,p_1,p_2\}$.
Then, for every vector, the width function of $K$ is strictly larger
than that of $P$, and furthermore, the lattice width of $P$ is equal
to $1$.
It follows that $\lwid > 1$.

If $K$ contains more than three integer points (and thus there is a
unimodular transformation of $K$ which contains $[0,1]^2$), the
assumptions of Lemma~\ref{bounds:for:quad} are fulfilled and the
assertion follows directly from Lemma~\ref{bounds:for:quad}.

Now assume that every edge of $K$ contains precisely one integer
point, i.e., $K$ is a triangle of type 3.
We define $K=Q=\conv \{q_0,q_1,q_2\}$ with $q_0,q_1,q_2$ and $Q$ given
as in Lemma~\ref{FreeCond}, and also borrow the other notations of
Lemma~\ref{FreeCond}.
Without loss of generality we assume that $x_0 \le x_1 \le x_2$ and
$x_i+x_j > 1$ for all $0 \le i < j \le 2$.
Let $f(x):=x_0x_1x_2+(1-x_0)(1-x_1)(1-x_2)$.
The upper bound for $\lwid$ follows from
Lemma~\ref{FreeCond}.

\emph{Case~1:} $x_0 \ge \frac{1}{2}$.
If $1 < \lwid < 2$, then in view of Lemma~\ref{tri-wid-analyticially}
we obtain 
$$
  A = \frac{1}{2 f(x)} = \frac{x_0}{f(x)} \cdot \frac{1}{2 x_0} = \frac{w}{2x_0} \le w < \frac{w^2}{2 (w-1)}.
$$
Assume now that $\lwid \ge 2$.
Then, taking into account Lemmas~\ref{FreeCond} and
\ref{tri-wid-analyticially} we obtain 
\begin{align*}
  A & \le \max \setcond{\frac{1}{2 f(x)}}{\frac{1}{2} \le x_0 \le
  x_1 \le x_2 < 1, \ x_0 = \lwid f(x)} \\
  & = \frac{1}{2} \left( \min \setcond{f(x)}{\frac{1}{2} \le x_0 \le
  x_1 \le x_2 <1, \ x_0 = \lwid f(x)} \right)^{-1}.
\end{align*}
Furthermore, using \eqref{detA lower}, we obtain
\begin{align*}
  &~\min \setcond{f(x)}{\frac{1}{2} \le x_0 \le x_1 \le x_2 < 1, \ x_0
    = \lwid f(x)} \\
  = &~\frac{1}{\lwid}\min \setcond{x_0}{\frac{1}{2} \le x_0 \le x_1 \le
    x_2 <1, \ x_0 = \lwid f(x)} \\
  \ge &~\frac{1}{\lwid} \min \setcond{x_0}{\frac{1}{2} \le x_0 < 1, \ x_0
    \ge \lwid (1- 3x_0+3 x_0^2)}.
\end{align*}
The minimum in the previous expression is attained for $x_0$ equal to
the smaller root of the equation $t = \lwid (1 - 3 t+ 3 t^2)$ since
this root, which is equal to
\[ \frac{3 \lwid + 1 - \sqrt{1+6 \lwid - 3 \lwid^2}}{6 \lwid }, \]
lies in  $[\frac{1}{2},1)$.
The characterization of the equality case follows directly by
analyzing the equality cases in the above estimates.
Thus, equality holds for $x_0 = x_1 = x_2$.

\emph{Case~2:} $x_0 < \frac{1}{2}.$ We have 
\begin{align*}
  I:=&~\inf \setcond{f(x)}{0 < x_0 < \frac{1}{2}, \, \frac{1}{2} <
  x_1 \le x_2 < 1, \, x_0+x_1 > 1, \, x_0+x_2 > 1, \, x_0 = \lwid
  f(x)} \\
  = &~\frac{1}{\lwid} \inf \setcond{x_0}{0 < x_0 < \frac{1}{2}, \,
  \frac{1}{2} < x_1 \le x_2 < 1, \, x_0+x_1 > 1, \, x_0+x_2 > 1, \,
  x_0 = \lwid f(x)} \\
  \ge &~\frac{1}{\lwid} \inf \setcond{x_0}{0 < x_0 < \frac{1}{2}, \,
  \frac{1}{2} < x_1 \le x_2 < 1, \, x_0+x_1 > 1, \, x_0+x_2 > 1, \,
  x_0 \ge \lwid f(x)}. 
\end{align*} 
Furthermore,
\begin{align*}
        f(x)  = &~(x_0 + x_1 - 1) x_2 + (1-x_0) (1-x_1) \\ 
               \ge &~(x_0 + x_1-1) x_1 + (1-x_0)(1-x_1) \\
                = &~x_0 (2 x_1 - 1) + (1-x_1)^2 \\
                = &~(1-x_0) x_0 + (x_0+x_1-1)^2 \\
                > &~(1-x_0) x_0. 
\end{align*}
Hence 
\begin{align*}
  I \ge &~\frac{1}{\lwid} \inf \setcond{x_0}{0 < x_0 < \frac{1}{2}, \,
  x_0 > w x_0 (1-x_0)} \\
  = &~\frac{1}{\lwid} \inf \setcond{x_0}{1- \frac{1}{\lwid} < x_0 <
  \frac{1}{2} } \\
         = &~\frac{1}{\lwid} \left( 1 - \frac{1}{\lwid}\right). 
\end{align*}
It follows that 
$$
        A = \frac{1}{2 f(x)} \le \frac{1}{2 I} \le \frac{\lwid^2}{2 (\lwid-1)}.
$$
The equality case is characterized in a straightforward way.
\end{proof}

We notice that Theorem~\ref{hurkens-thm} is a consequence of
Propositions~\ref{list-of-two-dim-max}, \ref{completing} and
Lemmas~\ref{bounds:for:quad}, \ref{bounds:for:tri} (in fact, the main
steps of the proof from \cite{Hurkens90} were incorporated in
Lemmas~\ref{FreeCond}, \ref{bounds:for:quad} and
\ref{bounds:for:tri}). 

\begin{proof}[Proof of Theorem~\ref{AreaLWidth}]
Let us show \eqref{InfBound}--\eqref{AreaHard}.
In view of Lemmas~\ref{bounds:for:quad}, \ref{bounds:for:tri} bounds
\eqref{InfBound}--\eqref{AreaHard} hold for maximal lattice-free
sets.
Let $K \in \KK^2$ be an arbitrary lattice-free set.
By Proposition~\ref{completing}, there exists a maximal lattice-free
set $H \in \KK^2$ such that $K \subseteq H$.
We have $\lwid(K) \leq \lwid(H)$ and $A(K) \leq
A(H)$.
For $0 < \lwid \le 1 + \frac{2}{\sqrt{3}}$ define $F(\lwid)$ to be the
upper bounds in \eqref{InfBound}--\eqref{AreaHard}.
Note that $F(\lwid)$ is monotonically non-increasing.
Thus, it follows $A(K) \le A(H) \le F(\lwid(H)) \le
F(\lwid(K))$.
The equality $A(K)=F(\lwid(K))$ implies $K=H$ and, thus, the
characterizations of the equality cases for
\eqref{InfBound}--\eqref{AreaHard} follow from the characterizations
of the equality cases in Lemmas~\ref{bounds:for:quad},
\ref{bounds:for:tri}.
The bound \eqref{AreaBelow} and Part~\ref{AreaBelowEq} follow directly
from Theorem~\ref{makaitoth-thm}.
\end{proof}

\begin{proof}[Proof of Corollary~\ref{Upper bounds for area in terms of mu's}]
We have $\mu_1(K) \wid(K) = 1$.
Furthermore, an appropriate translate of $\mu_2(K) \cdot K$ is lattice-free.
Thus, we apply the upper bounds of Theorem~\ref{AreaLWidth} to the
body $\mu_2(K) \cdot K$ and then express the lattice width of
$\mu_2(K) \cdot K$ as
$\frac{\mu_2(K)}{\mu_1(K)}.$
The sharpness of the bounds follows from the characterizations of the
equality cases in Theorem~\ref{AreaLWidth}.
\end{proof}

%%%%%%%%%%%%%%%%%%%%%%%%%%%%%%%%%%%%%%%%%%%%%%%%%%%%%%%%%%%%%%%%%%%%%%
%%%%%%%%%%%%%%%%%%%%%%%%%%%%%%%%%%%%%%%%%%%%%%%%%%%%%%%%%%%%%%%%%%%%%%

\section{Proofs for centrally symmetric bodies} \label{centr.sym.case}

In this section we prove Theorems \ref{CentrSym} and
\ref{cov min c-sym planar}.
Inequalities \eqref{SymLWidthUpper} have already been stated in
\cite[Theorem~2.13]{KannanLovasz88}, but the proof of Kannan and
Lov\'asz does not seem to show this result.
Therefore, we first show \eqref{SymLWidthUpper} by proving Theorem
\ref{cov min c-sym planar}. Afterwards, we use \eqref{SymLWidthUpper}
to show \eqref{SymInf}--\eqref{SymAreaLower}.

\begin{figure}[ht]
  \unitlength=1mm
  \centering
  \subfigure[\label{3-hexagons-1st-reduction} Reduction of $R$ (thin solid) to a
    parallelogram $R'$ (dashed); $R$ and $R'$ are circumscribed about $Q$ (bold solid); $Q$ is circumscribed about $P$ (dotted)]{
    \begin{picture}(70,70)
      \put(0,3){\includegraphics[width=70mm]{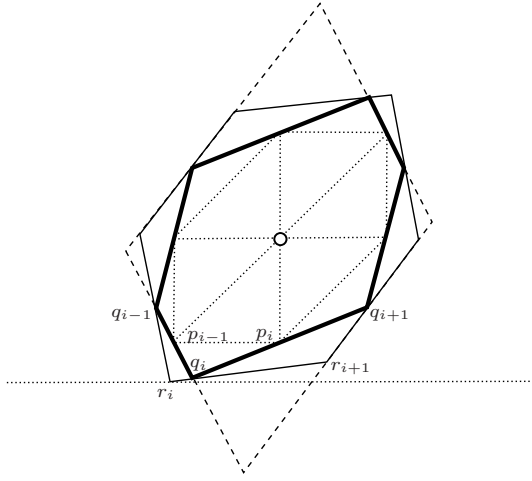}}
      \put(20,13.4){\scriptsize $r_i$}
      \put(42.8,17){\scriptsize $r_{i+1}$}
      \put(24.3,17.3){\scriptsize $q_i$}
      \put(48,24){\scriptsize $q_{i+1}$}
      \put(33,21.3){\scriptsize $p_i$}
      \put(24,21.1){\scriptsize $p_{i-1}$}
      \put(14,24){\scriptsize $q_{i-1}$}
      %\graphpaper[5](0,0)(70,70)
    \end{picture}}
  \qquad \qquad
  \subfigure[\label{3-hexagons-2nd-reduction} Reduction of $Q$ (bold solid) and $R'$ (thin solid) to 
    parallelograms $Q'$ (bold dashed) and $R''$ (thin dashed), respectively]{
    \begin{picture}(70,70)
      \put(15,1){\includegraphics[width=40mm]{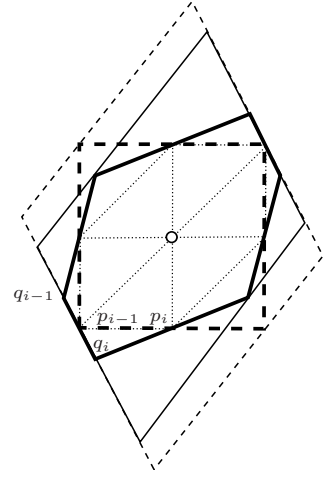}}
      \put(14,24){\scriptsize $q_{i-1}$}
      \put(24.5,17.3){\scriptsize $q_i$}
      \put(25,21){\scriptsize $p_{i-1}$}
      \put(32,21){\scriptsize $p_i$}
      %\graphpaper[5](0,0)(70,70)
    \end{picture}
  }
  \caption{ The hexagons $Q$ and $R$ degenerate into
    parallelograms\label{3-hexagons-reduction}}
\end{figure}

\begin{proof}[Proof of Theorem~\ref{cov min c-sym planar}]
We want to show that $\mu_2(K) \le 2\mu_1(K)$ for every centrally
symmetric $K \in \KK^2$.
For convenience we define $\mu_i := \mu_i(K)$ for $i = 1,2$ and
$\lwid := \lwid(K)$. We further assume that $K$ has $o$ as its center
of symmetry.
Since both $\mu_1$ and $\mu_2$ are homogeneous of degree $-1$, it
suffices to consider the case $\mu_2 = 1$.
Using $\mu_1 \lwid=1$ it remains to prove that $\lwid \le 2$.

In this proof, when we say \term{hexagon} we refer to a centrally
symmetric convex hexagon which could also be degenerated so as to be
a parallelogram.
A hexagon will be given by a sequence of consecutive vertices indexed
modulo $6$.
Because of degeneracy we allow the ``vertices'' to be collinear or to
coincide.

Since $\mu_2 = 1$, the translates $K + z$ with $z \in \integer^2$
cover $\real^2$.
By \cite[Lemma~3, p.~246]{GruLek87} there exists a hexagon $Q := \conv
\setcond{q_i}{i=0,\ldots,5}$ centered at $o$ such that the vertices of
$Q$ lie on the boundary of $K$ and the translates $Q + z$ with $z \in
\integer^2$ tile $\real^2$.
We introduce a second hexagon $R:= \conv \setcond{r_i}{i=0,\ldots,5}$,
centered at $o$ and containing $Q$, such that for every $i$ the line
$\aff \{r_i,r_{i+1}\}$ supports $K$ at $q_i$. This implies $\mu_2(R) =
1$ and from $K \subseteq R$ it follows that $\lwid\le \lwid(R)$. 
Thus, it suffices to bound $\lwid(R)$ from above. 

Let $P := \conv \setcond{p_i}{i=0,\ldots,5}$ be the hexagon with $p_i
:= \frac{1}{2}(q_i+q_{i+1})$ for every $i$.
From the choice of $Q$ it follows that $2 p_i, 2 p_{i+1}$ is a basis
of $\integer^2$ for every $i$.
Thus, we constructed hexagons $P \subseteq Q \subseteq R$ such that
every vertex of $P$ is the central point of an edge of $Q$, and each
edge of $R$ contains a vertex of $Q$ (see
Fig.~\ref{3-hexagons-1st-reduction}). 
We show that for such a configuration there exists an edge of $P$ such
that its normal, say $u$, satisfies $w(R,u) \le 2 w(P,u)$.
In terms of lattices this means that $\lwid(R) \le 2$, which then
finishes the proof.

By $u_i$ we shall denote the normal vector to the edge $\conv
\{p_i,p_{i+1}\}$ of $P$.
Let us take any vertex $r_i$ of $R$ such that some line through $r_i$
which is parallel to an edge of $P$ supports $R$. This supporting line
is then parallel to either $\aff \{p_{i-1},p_i\}$ or $\aff
\{p_{i-2},p_{i-1}\}$.
By symmetry we may assume that the line through $r_i$ parallel to
$\aff \{p_{i-1},p_i\}$ supports $R$.
Now consider the parallelogram $R'$ which is the intersection of the
four halfplanes defined by the lines $\pm \aff \{q_{i-1},q_i\}$ and
$\pm \aff \{r_{i+1},r_{i+2}\}$ containing $o$.
It can be seen that $w(R,u_j) \le w(R',u_j)$ for every $j$.
Thus, we can replace $R$ by $R'$ and it suffices to show $\lwid(R') \le 2$. 

The edge $\conv \{q_{i-1},q_i\}$ of $Q$ lies on the boundary of $R'$.
We modify $Q$ and $R'$ simultaneously to obtain $Q'$ and $R''$,
respectively, in the following way.
Define $Q' := \conv \{ \pm p_{i-1}, \pm (p_i+p_{i+1})\}$ and define
$R''$ as the parallelogram circumscribing $Q'$ with edges parallel to
the edges of $R'$ (see Fig.~\ref{3-hexagons-2nd-reduction}).
Clearly, $R' \subseteq R''$ and by this $w(R',u_j) \le w(R'',u_j)$ for
every $j$.

With our modifications we have reduced the problem on hexagons to a
problem on quadrilaterals, i.e., where both $Q$ and $R$ are
quadrilaterals.
Now, the assertion follows using the same arguments as in the proof of
Lemma~\ref{bounds:for:quad}. There is a unimodular transformation
which maps $Q$ to $[0,1]^2$ and $R$ to a quadrilateral $K = \conv
\{q_1,q_2,q_3,q_4\}$ as in Lemma~\ref{bounds:for:quad} (see also
Fig.~\ref{QuadFig}), where we have showed that $w(K,u) \le 2$ for
$u = (1,0)$ or $u = (0,1)$.
\end{proof}

 Remember that every $o$-symmetric $K \in \KK^2$ can be associated with a norm 
$$\|u\|_K = \min \setcond{\alpha\ge 0}{u \in \alpha K}.$$ It is well known that $\|u\|_K = h(K^\ast,u)$.

\begin{lemma}\label{mu_2 for alpha-diamond}
 Let $0 \le \alpha < 1$ and let $K_\alpha:=\conv \{\pm (1,\alpha),\pm (0,1)\}$. Then $\mu_2(K_\alpha) = \frac{1}{2} \max \{1+\alpha,2-\alpha\}$.
\end{lemma}
\begin{proof}
The unimodular transformation given by the matrix $\begin{pmatrix} 1 & 0 \\ 1 & -1 \end{pmatrix}$ maps $K_\alpha$ onto $K_{1-\alpha}$. Thus, it suffices to consider the case $0 \le \alpha \le \frac{1}{2}$. For the sake of brevity we write $K:=K_\alpha$.
Direct computations show that $K^\ast = \conv \{\pm (-\alpha+1,1),\pm (-\alpha-1,1)\} = \conv \{(-\alpha,1),(\alpha,-1)\} +\conv \{(1,0),(-1,0)\}$.
Hence 
$$f(u):=h(K^\ast,u) = |u_2-\alpha u_1|+|u_1|,$$
for $u=(u_1,u_2) \in \real^2$.
We have 
\begin{align}
 \mu_2(K) 
        & = \min \setcond{\mu\ge 0}{\mu K + \integer^2=\real^2} \nonumber\\
        & = \min \setcond{\mu \ge 0 }{\forall x \in \real^2 \, \exists z\in \integer^2 \ \mbox{such that} \ \|x-z\|_K \le \mu} \nonumber\\
        & = \max_{x \in \real^2} \min_{z\in \integer^2} \|x-z\|_K \nonumber\\
        & = \max_{x \in \real^2} \min_{z\in \integer^2} h(K^\ast,x-z) \nonumber\\
        & = \max_{x \in \real^2} \min_{z\in \integer^2} f(x-z).\label{mu maximin}
\end{align}
For $t \in \real$ let $\nint{t}$ denote the \term{nearest integer function}. If $s,t \in \real$, we introduce the distance between $s$ and $t$ modulo $1$ by $d(s,t):=d(s-t,0),$ where $d(s,0):=|s-\nint{s}|$. We choose $x=(x_1,x_2) \in \real^2$ and consider
$$
        f(x-z)=|x_2-z_2-\alpha(x_1-z_1)|+|x_1-z_1|.
$$
with $z=(z_1,z_2)$ varying in $\integer^2$. For $z_1=\nint{x_1}$ and $z_2=\nint{x_2-\alpha(x_1-z_1)}$ we have $f(x-z) \le 1$. Furthermore, if $z_1 \not\in \left\{ \floor{x_1}, \ceil{x_1} \right\}$, then $f(x-z) \ge |x_1-z_1| \ge 1$. Thus, computing \eqref{mu maximin} we may assume that $z_1 \in \left\{\floor{x_1}, \ceil{x_1} \right\}$. If $x_1 \in \integer$ we may set $z_1=x_1$ and $z_2=\nint{x_2}$ obtaining $f(x-z) \le \frac{1}{2}$. Otherwise $x_1 \in \real \setminus \integer$ and we introduce $\beta:=x_1-\floor{x_1}$ satisfying $0 < \beta < 1$. Since we can assume that $z_1 \in \left\{ \floor{x_1}, \ceil{x_1} \right\}$ we have 
\begin{align*}
 \min_{z\in \integer^2} f(x-z) 
        &= \min_{z_2 \in \integer} \min \left\{|x_2-z_2-\alpha \beta|+\beta,|x_2-z_2+\alpha(1-\beta)|+1 - \beta\right\} \\
        & = \min \left\{d(x_2-\alpha\beta,0)+\beta,d(x_2+\alpha(1-\beta),0)+1-\beta \right\}.
\end{align*}
Thus, the value in \eqref{mu maximin} is the maximum of $\frac{1}{2}$  and the value 
\begin{align*}
\max_{\overtwocond{x \in \real^2}{x_1 \not\in \integer}}\min_{z \in \integer^2} f(x-z) 
        & = \max_{\overtwocond{0 < \beta< 1}{x_2 \in\real}} \min \left\{d(x_2-\alpha\beta,0)+\beta,d(x_2+\alpha(1-\beta),0)+1-\beta \right\} \\
                        & = \max_{\overtwocond{0 < \beta< 1}{x_2 \in\real}} \min \left\{d(x_2+\alpha(1-\beta),\alpha)+\beta,d(x_2+\alpha(1-\beta),0)+1-\beta \right\} \\
        &= \max_{\overtwocond{0 < \beta< 1}{x_2 \in\real}} \min \left\{d(x_2,\alpha)+\beta,d(x_2,0)+1-\beta \right\} \\
        & \le \max_{\beta,x_2 \in \real} \min \left\{d(x_2,\alpha)+\beta,d(x_2,0)+1-\beta \right\}.
\end{align*}
If $d(x_2,\alpha)+\beta$ and $d(x_2,0)+1-\beta$ differ, then sligthly perturbing $\beta$ the minimum of these two values becomes larger. Hence the latter maximum is attained for the $\beta$ for which $d(x_2,\alpha)+\beta$ and $d(x_2,0)+1-\beta$ coincide. In this case $\beta= \frac{1}{2}(d(x_2,0)-d(x_2,\alpha)+1).$ Thus
\begin{align*}
 \max_{\beta,x_2 \in \real} \min \left\{d(x_2,\alpha)+\beta,d(x_2,0)+1-\beta \right\} & = \max_{x_2 \in \real} \frac{1}{2}(d(x_2,0)+d(x_2,\alpha)+1) \\ & = \frac{1}{2} \Big(1+ \max_{x_2 \in \real} \big(d(x_2,0)+d(x_2,\alpha)\big) \Big) \\
        & = \frac{1}{2}(2-\alpha).
\end{align*}

Since the latter is at least $\frac{1}{2}$ we have shown $\mu_2(K) \le \frac{1}{2} (2-\alpha).$ It remains to show that the above is attained with equality. Employing the above derivations we can see that the equality is attained for $x_1=\frac{1+\alpha}{2}$ and $x_2=\frac{1}{2}$.
\end{proof}

\begin{proof}[Proof of Theorem~\ref{CentrSym}]
Since \eqref{SymLWidthUpper} is established the bounds \eqref{SymInf}
and \eqref{SymAreaUpper} together with
Parts~\ref{eq case for inf bound, sym} and
\ref{eq case for area bound} follow directly from
Theorem~\ref{AreaLWidth}. 

Now let us show \eqref{SymAreaLower} and Part
\ref{eq case lower bound sym}.
By \eqref{lattice-width-over-lambda}, $o$ is the only interior integer
point in $\lwid \cdot (DK)^\ast$, where $\lwid:=\lwid(K)$.
Thus, by Minkowski's first theorem, $\area(\lwid \cdot (DK)^\ast) \le 4$.
Using the above fact and Mahler's inequality we obtain
\begin{align*}
  \area(K) = \frac{\area(DK)}{4} = \frac{\area(DK) \area((DK)^\ast)}{4
  \area ((DK)^\ast)} \ge \frac{2}{\area((DK)^\ast)} \ge
  \frac{\lwid^2}{2}.
\end{align*}

It remains to characterize the case where this inequality is tight.
In view of Mahler's inequality and Minkowski's first theorem (the
parts in these statements which give information on the equality
cases), the equality $\area(K) = \frac{\lwid^2}{2}$ implies that
$(DK)^\ast$ is a parallelogram and that the sets $\frac{\lwid}{2}
(DK)^\ast + z$ with $z \in \integer^2$ tile $\real^2$.
By Proposition~\ref{tiling-by-parallelograms} we deduce that there
exists a linear unimodular transformation $T$ and a parameter $0 \le \alpha < 1$ such that
\[ T\left( \frac{\lwid}{2} (DK)^\ast \right) = \frac{1}{2} \conv
   \{\pm(-\alpha-1,1),\pm(-\alpha+1,1)\}. 
\]
Taking duals and slightly modifying the left hand
side of the equality we arrive at
\[ \frac{2}{\lwid} (T^{-1})^\ast (DK) = 2 \left( \conv \{\pm
   (-\alpha-1,1),\pm (-\alpha+1,1) \}\right)^\ast.
\]
Direct computations yield
\[ \left( \conv \{\pm (-\alpha-1,1),\pm (-\alpha+1,1) \}\right)^\ast =
   \conv \{ \pm (1,\alpha) , \pm (0,1) \}.
\]
Clearly, the transformation $(T^{-1})^\ast$ is unimodular.
Summarizing we see that, up to unimodular transformations, we have 
\[ \frac{1}{2} DK = \frac{\lwid}{2} \conv \{ \pm(1,\alpha),  \pm (0,1) \} \]
with $0 \le \alpha < 1$.
Note that $\frac{1}{2} DK$ is a translate of $K$. An appropriate translate of $\frac{1}{2} DK$ is lattice-free if and only if $\mu_2(\frac{1}{2} DK) \ge 1.$ In view of Lemma~\ref{mu_2 for alpha-diamond} we get $\mu_2(\frac{1}{2} DK) = \frac{2}{\lwid} \cdot \frac{1}{2} \max \{1+\alpha,2-\alpha\}$. The above observations yield Part~\ref{eq case lower bound sym}.
\end{proof}

\begin{proof}[Proof of Corollary~\ref{mu.sym}]
        The proof is analogous to the proof of Corollary~\ref{Upper bounds for area in terms of mu's}.
\end{proof}

\subsection*{Acknowledgements}

We thank the anonymous referee, E.~Makai, Jr.~and
C.~A.~J.~Hurkens for useful suggestions and pointers to the
literature.

%%%%%%%%%%%%%%%%%%%%%%%%%%%%%%%%%%%%%%%%%%%%%%%%%%%%%%%%%%%%%%%%%%%%%%
%%%%%%%%%%%%%%%%%%%%%%%%%%%%%%%%%%%%%%%%%%%%%%%%%%%%%%%%%%%%%%%%%%%%%%

%\bibliographystyle{amsplain}
%\bibliography{latconv}
\providecommand{\bysame}{\leavevmode\hbox to3em{\hrulefill}\thinspace}
\providecommand{\MR}{\relax\ifhmode\unskip\space\fi MR }
\providecommand{\MRhref}[2]{%
  \href{http://www.ams.org/mathscinet-getitem?mr=#1}{#2}
}
\providecommand{\href}[2]{#2}

\begin{tabular}{l}
  \textsc{Gennadiy Averkov} \\
  \textsc{Institute of Mathematical Optimization} \\
  \textsc{Faculty of Mathematics} \\
  \textsc{University of Magdeburg} \\
  \textsc{Universit\"atsplatz 2, 39106 Magdeburg} \\
  \textsc{Germany}        \\
  \emph{e-mail}: \texttt{averkov@ovgu.de} \\
  \emph{web}: \texttt{http://fma2.math.uni-magdeburg.de/$\sim$averkov}
  \\ \\
\end{tabular}

\begin{tabular}{l}
  \textsc{Christian Wagner} \\
  \textsc{Institute for Operations Research} \\
  %\textsc{Faculty of Mathematics} \\
  \textsc{ETH Z\"urich} \\
  \textsc{R\"amistrasse 101, 8092 Z\"urich} \\
  \textsc{Switzerland}        \\
  \emph{e-mail}: \texttt{christian.wagner@ifor.math.ethz.ch} \\
  \emph{web}: \texttt{http://www.ifor.math.ethz.ch/staff/chwagner}
\end{tabular}

\end{document}